\documentclass[11pt,a4paper]{article}
\usepackage[T1]{fontenc}
\usepackage[utf8]{inputenc}
\usepackage[hmargin={26mm,26mm},vmargin={30mm,35mm}]{geometry}
\usepackage[style=numeric-comp,maxbibnames=100,bibencoding=utf8,giveninits=true,backend=bibtex]{biblatex}
\usepackage{amsmath,amsfonts,amssymb,amsthm}
\usepackage{authblk}
\usepackage{array}
\usepackage{booktabs}
\usepackage{comment}
\usepackage[tablesonly]{endfloat}
\usepackage{enumitem}
\usepackage[colorlinks=true, allcolors=blue]{hyperref}
\usepackage{listings}
\usepackage{mathtools}
\usepackage[expansion=false]{microtype}
\usepackage{newtxtext,newtxmath}
\usepackage{polytopal-dsl-listings}
\usepackage{pgfplots}
\usepackage{subcaption}
\usepackage{tabularx}
\usepackage{xcolor}
\usepackage{tikz}
\usetikzlibrary{arrows.meta,backgrounds,fit,positioning}
\usepackage{xspace}

\pgfplotsset{compat=1.18}

\bibliography{poly_rust-paper}
\AtBeginBibliography{\small}

\newcommand{\email}[1]{\href{mailto:#1}{#1}}

\usepackage[normalem]{ulem}
\newcounter{corr}
\definecolor{violet}{rgb}{0.580,0.,0.827}
\newcommand{\corr}[3]{\typeout{Warning : a correction remains in page \thepage}
  \stepcounter{corr}
 	      {\color{blue}\ifmmode\text{\,\sout{\ensuremath{#1}}\,}\else\sout{#1}\fi}
              {\color{red}#2}
              {\color{violet} #3}
}

\lstdefinestyle{GenericStyle}{
  basicstyle=\ttfamily\small,
  language={},
  keywordstyle={},
  commentstyle={},
  stringstyle={}
}

\newcommand{\Real}{\mathbb{R}}

\newcommand{\dirichlet}{\mathrm{D}}
\newcommand{\neumann}{\mathrm{N}}

\newcommand{\Mh}{\mathcal{M}_h}
\newcommand{\Th}{\mathcal{T}_h}
\newcommand{\Eh}{\mathcal{E}_h}
\newcommand{\ET}{\mathcal{E}_T}

\newcommand{\Poly}[1]{\mathcal{P}^{#1}}

\newcommand{\norm}[2]{\|#2\|_{#1}}

\theoremstyle{plain}
\theoremstyle{remark}

\definecolor{grammarbackground}{RGB}{247,248,250}
\definecolor{grammarframe}{RGB}{205,210,218}
\definecolor{terminalcolour}{RGB}{35,78,140}

\hypersetup{
  colorlinks=true,
  linkcolor=terminalcolour,
  urlcolor=terminalcolour,
  pdftitle={Polytopal DSL Syntax Reference}
}

\newcommand{\term}[1]{\texttt{\color{terminalcolour}#1}}

\newcommand{\polyrust}{\texttt{poly\_rust}\xspace}
\newcommand{\Rust}{\texttt{Rust}\xspace}
\newcommand{\cargo}{\texttt{cargo}\xspace}
\newcommand{\Cpp}{\texttt{C++}\xspace}

\title{\polyrust: a domain-specific language for polytopal methods with a \Rust interpreter\footnote{arXiv preprint \href{https://arxiv.org/2608.09414}{2608.09414 [math.NA]}}}
\author{Daniele Antonio Di Pietro}
\affil{%
  Universit\'e de Montpellier, Montpellier 34090, France \\
  \email{daniele.di-pietro@umontpellier.fr} \\
  Corresponding author
}

\begin{document}

\maketitle

\begin{abstract}
  \noindent\textbf{Background:} Over the last two decades, considerable effort in the applied mathematics community has been devoted to the development of polytopal methods, i.e., numerical methods for partial differential equations that support meshes more general than those used by traditional finite element methods.
  Despite their flexibility and attractive mathematical properties, polytopal methods have a steep implementation curve and, to date, no easily accessible open-source software framework has been available for their systematic implementation.
  \smallskip\\
  \textbf{Methods:} In this article, we introduce \polyrust, a new software library implementing a domain-specific language (DSL) specifically designed for polytopal methods.
  \polyrust provides high-level constructs for expressing a wide range of such methods using a syntax that closely mirrors their mathematical formulation, thereby substantially lowering the implementation barrier for both newcomers and experienced users.
  \smallskip\\
  \textbf{Results:} The library facilitates the rapid implementation, testing, and comparison of polytopal discretizations, making these methods more accessible to researchers, applied scientists, engineers, and students.
  In our vision, this will foster their dissemination beyond the numerical analysis community and encourage their adoption by practitioners in a broad range of application fields.
  It will also contribute to reproducibility of results, as DSL files are text files that can be easily made available, stored, and exchanged.
  \smallskip\\
  \textbf{Conclusions:} \polyrust introduces, to our knowledge, the first DSL specifically tailored to polytopal methods.
  By combining mathematical expressiveness, ease of implementation, and an open-source framework, it provides a common platform for developing, sharing, reproducing, and comparing polytopal schemes, with the ambition of becoming a reference tool for the wider scientific community.
\end{abstract}


\section{Introduction}

Over the last twenty years, significant research efforts have been made to develop and analyze discretization methods for partial differential equations that support more general meshes than those on which traditional finite element methods are built.
A finite element is classically defined as a triplet $(T, V, \sigma)$, where $T$ is the geometric element, $V$ is a finite-dimensional space of functions over $T$, and $\sigma$ is a set of degrees of freedom (DOFs), i.e., a basis of the dual space of $V$.
The DOFs are selected so that appropriate global continuity properties result from their single-valuedness at internal vertices, edges, and faces.
Constructing $V$ and $\sigma$ so that the above requirements are fulfilled is generally possible only when $T$ has a simple shape, e.g., a triangle or a rectangle in dimension two.

In order to treat geometric elements $T$ of general shape, polytopal methods renounce conformity and often drop the space $V$ altogether, building the relevant discrete operators by direct manipulations of the DOFs (a non-computable space $V$ may still be used for the sake of presentation \cite{Beirao-da-Veiga.Brezzi.ea:13}).
This leads to major implementation differences with respect to finite element methods, as constructing these operators typically requires solving local variational problems on appropriately selected mesh entities.
A variety of polytopal methods exist, with differences that are both philosophical and implementative.
We refer to, e.g., \cite{Di-Pietro.Droniou:20} for a broad introduction to the subject and a mostly up-to-date bibliography.
\polyrust is mainly focused on methods that fall within the broad categories of Virtual Elements \cite{Beirao-da-Veiga.Brezzi.ea:13,Beirao-da-Veiga.Brezzi.ea:14,Beirao-da-Veiga.Brezzi.ea:23},
Hybrid High-Order (HHO) \cite{Di-Pietro.Ern.ea:14,Di-Pietro.Ern:15,Di-Pietro.Droniou:20},
or Discrete de Rham (DDR) methods \cite{Di-Pietro.Droniou.ea:20,Di-Pietro.Droniou:23*2,Bonaldi.Di-Pietro.ea:25}.
The basic tools it offers also enable the implementation of polytopal discontinuous Galerkin methods \cite{Di-Pietro.Ern:10,Bassi.Botti.ea:12,Antonietti.Giani.ea:13,Antonietti.Cangiani.ea:16}, but this point is less original and will not be treated in detail in the present work.
\smallskip

\begin{figure}
  \centering
  \includegraphics[width=0.8\textwidth]{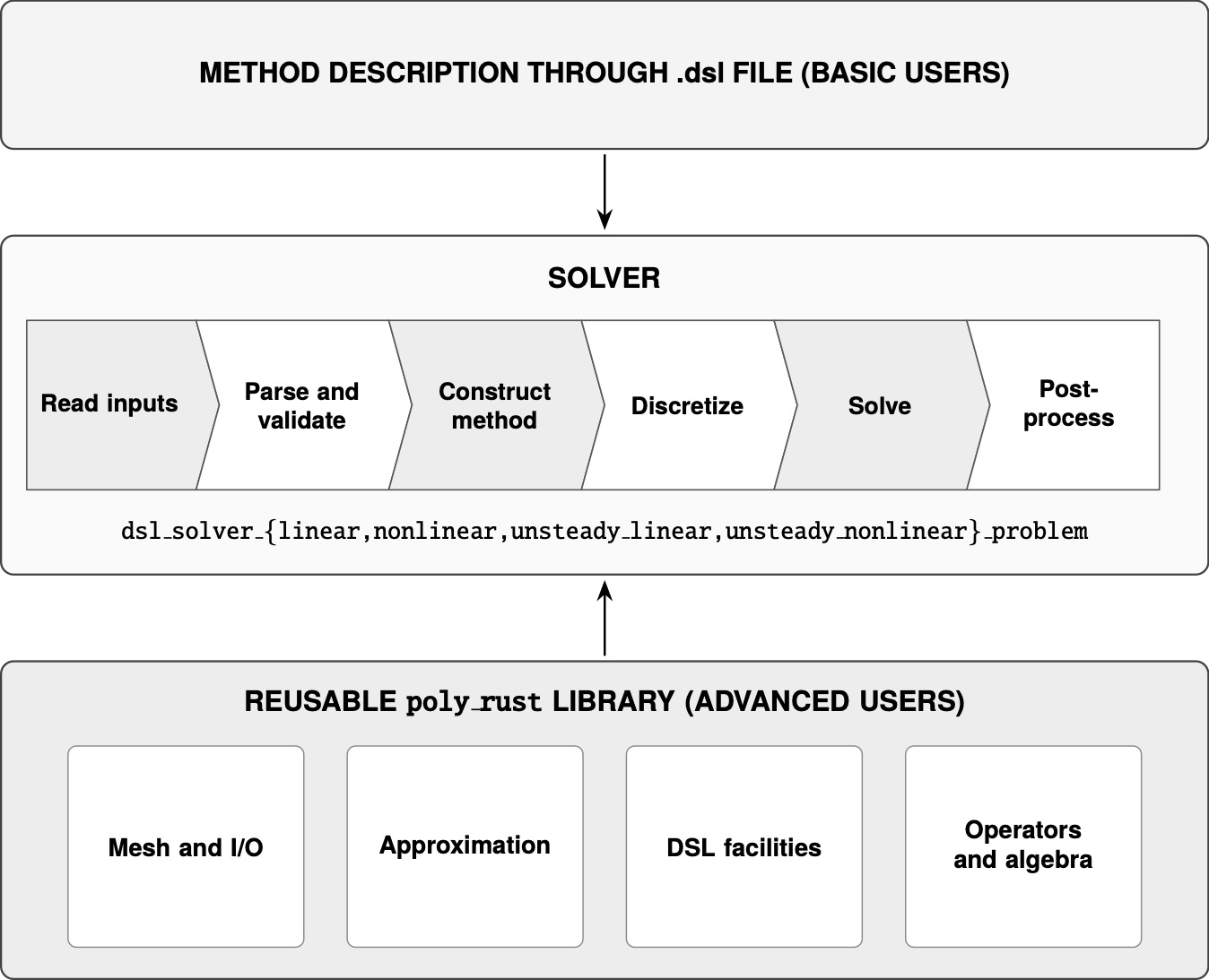}
  \caption{Layered architecture of \polyrust.}
  \label{fig:poly-rust-architecture}
\end{figure}
\polyrust provides a domain-specific language (DSL) of broad scope for the implementation of polytopal methods.
The adopted vision is fully discrete and DOF-centered, and the definitions of the relevant objects closely resemble their mathematical counterparts.
The architecture of \polyrust is depicted in Figure \ref{fig:poly-rust-architecture}.
An application is composed of two main ingredients:
\begin{itemize}
\item A \emph{DSL configuration file}, describing the discrete problem;

\item A \emph{solver}, which processes the DSL, builds and solves the discrete problem, and produces the requested output.
\end{itemize}
Solvers are built upon the reusable ingredients of \polyrust, which include basic I/O, a numerical core, DSL-related facilities, and the assembly and solution of sparse algebraic systems.
Basic users describe a method in a DSL file and run one of the supplied generic solvers.
A solver orchestrates method construction, assembly, solution, and post-processing using the reusable library facilities.
Advanced users can compose the same facilities into new generic or problem-specific solvers.
Our focus is, for the sake of simplicity, on two-dimensional problems.
Support for three-dimensional problems along with other generic solvers are available in the current version v0.13.1.

The language for solvers and library components is \Rust\footnote{\url{https://rust-lang.org}}.
\Rust is not (yet) common in scientific computing, but it has several advantages compared with more established languages.
Among the most important ones, we can mention:
the availability of comprehensive online documentation;
performance and memory safety by design;
powerful native types;
less compiler-intensive handling of generic types compared with, e.g., \Cpp templates;
a native package manager and build tool, \cargo.
The implementation relies on several crates available from \url{https://crates.io}
\texttt{ndarray} and \texttt{ndarray-linalg} for dense arrays and local linear algebra;
\texttt{quadraturerules} for quadrature rules on standard geometries;
\texttt{pest} to generate the parser;
\texttt{sprs} for sparse matrix storage and assembly;
\texttt{rayon} for shared-memory parallelism.
The dependencies of \polyrust external to \url{https://crates.io} are mainly optional and kept to a minimum in order to simplify compilation and deployment.

The purpose of this paper is not to provide an in-depth, comprehensive description of \polyrust, but rather to illustrate the main ideas and give basic users a glimpse of the underlying philosophy.
For this reason, we have preferred to introduce the syntax through meaningful examples rather than by a formal description.
Numerical examples covering diverse problems and methods complete the exposition.
\smallskip

Software support for the discretization of partial differential equations is extensive.
General-purpose finite element software spans several levels of abstraction.
\texttt{deal.II} is an object-oriented C++ library providing comprehensive building blocks for the development of adaptive finite element codes \cite{Bangerth.Hartmann.ea:07}.
\texttt{FreeFEM++} is an integrated environment with its own high-level language for defining meshes, finite element spaces, variational problems, and solution algorithms \cite{Hecht:12}, %
while \texttt{Feel++} provides an embedded domain-specific language in \Cpp for expressing function spaces, differential operators, integrals, and Galerkin formulations in a syntax close to their mathematical description \cite{Prudhomme.Chabannes.ea:12}.
UFL/FEniCS similarly represents weak forms through a domain-specific language embedded in Python and relies on automated code generation, %
whereas Gridap provides a high-level Julia API for writing weak formulations in near-mathematical notation, using Julia's just-in-time compilation rather than a separate form compiler \cite{Alns.Logg.ea:14,Badia.Verdugo:20}.

More closely related to the present work in terms of targeted methods,
\texttt{Vem++} \cite{Dassi:26} and PolyDiM \cite{Berrone.Borio.ea:26} focus on Virtual Elements,
\texttt{Disk++} \cite{Cicuttin.Di-Pietro.ea:18*1,Cicuttin.Ern.ea:21} and \texttt{ParaSkel++} \cite{Beaude.Lemaire:21} on HHO methods,
whereas HArDCore (\url{https://github.com/jdroniou/HArDCore}) offers reusable components for HHO and DDR methods on polytopal meshes.
We also cite here \texttt{Vulpes} (\url{https://vulpeslib.github.io}), mainly focusing on polytopal discontinuous Galerkin methods.
To the best of our knowledge, however, none combines general polytopal meshes with an external DSL covering all the required ingredients to express local reconstruction operators as well as discrete formulations that use them.
This is the main distinguishing feature of \polyrust and offers a great flexibility in the fast prototyping of new methods and the development of solvers for a large variety of partial differential equation problems.
\smallskip

The rest of this work is organized as follows.
In Section \ref{sec:overview} we provide an overview of the language, using the HHO method for the Poisson problem as a driving example.
In Section \ref{sec:solvers} we briefly review the problem-agnostic solvers provided for linear and nonlinear problems, along with other facilities useful for prototyping numerical methods.
Some advanced features are then discussed through examples in Section \ref{sec:advanced.features}, and Section \ref{sec:conclusion} concludes the paper.


\section{Overview of the language}\label{sec:overview}

In this section we discuss the DSL configuration file for a model problem, while the next section briefly presents the included basic solvers from a user's perspective.
The presentation follows the dependencies between the main language constructs listed in the top block of Figure \ref{fig:poly-rust-architecture}:
parameters and functions provide the problem data; discrete spaces and interpolants define the DOFs and their meaning; operators, forms, and boundary conditions define the discrete equations; and problem blocks select the ingredients to assemble.
Finally, post-processing is parametrized using functionals and export directives.

Denote by $\Omega \subset \Real^2$ a polygonal domain with boundary $\Gamma$ partitioned as follows:
\begin{equation}\label{eq:Gamma}
  \Gamma = \Gamma^\dirichlet \sqcup \Gamma^\neumann,
\end{equation}
where $\Gamma^\dirichlet$ is a relatively open subset of $\Gamma$ with positive one-dimensional Hausdorff measure.
Let $f : \Omega \to \Real$ be a source function,
$g : \Gamma^\dirichlet \to \Real$ the Dirichlet datum on $\Gamma^\dirichlet$,
and $\phi : \Gamma^\neumann \to \Real$ the Neumann datum on $\Gamma^\neumann$.
We will focus on the Hybrid High-Order (HHO) discretization of the Poisson problem:
Find $u : \Omega \to \Real$ such that
\begin{equation}\label{eq:continuous}
  \begin{alignedat}{2}
    - \Delta u &= f &\qquad& \text{in $\Omega$},
    \\
    u &= g &\qquad& \text{on $\Gamma^\dirichlet$},
    \\
    -\nabla u \cdot n &= \phi &\qquad& \text{on $\Gamma^\neumann$},
  \end{alignedat}
\end{equation}
where $n$ denotes the outward unit normal to $\Gamma$.
In what follows, we denote by $\Mh \coloneqq \Th \cup \Eh$ a polygonal mesh of $\Omega$ in the sense of \cite[Definition 1.4]{Di-Pietro.Droniou:20}, with $\Th$ collecting mesh elements and $\Eh$ mesh edges.
We assume that the trace of $\Mh$ on $\Gamma$ is compatible with the partition \eqref{eq:Gamma}.
In order to showcase a full pipeline including error computations, we consider a manufactured analytical solution, but this is of course not necessary.

\subsection{Method}

All declarations describing the numerical method for an application (and, in particular, all the blocks discussed in the rest of this section) are contained in an overarching \lstinline{method} block:
\begin{lstlisting}
method hho_poisson {
  // Parameters and functions
  // Discrete spaces
  // Interpolants
  // Operators
  // Forms
  // Norm functionals
  // Boundary conditions
  // Problem description
}
\end{lstlisting}
The method declaration can contain conditions on the polynomial degree \lstinline{k} to account, e.g., for the fact that a method is defined/convergent only for certain values.

\subsection{Parameters and functions}

The DSL permits the definition of parameters and functions, as exemplified in the following listing.
\begin{lstlisting}[caption={Parameters and functions}, label=lst:parameters.functions]
parameter pi = 3.141592653589793238462643383279502884197169399375105820974944592307816406286198

function exact_solution(vector X) -> scalar = sin(pi * X[0]) * sin(pi * X[1])
function exact_solution_gradient(vector X) -> vector
  = vector(
           pi * cos(pi * X[0]) * sin(pi * X[1]),
           pi * sin(pi * X[0]) * cos(pi * X[1])
           )
function exact_normal_flux(vector X) -> scalar on edge E of element T
  = - exact_solution_gradient(X) dot normal
           
function source(vector X) -> scalar = 2.0 * pow(pi, 2.0) * sin(pi * X[0]) * sin(pi * X[1])
\end{lstlisting}
Both parameters and functions can be scalar-, vector-, or matrix-valued, and they can also depend on the geometric context, as does the \lstinline{exact_normal_flux} function, which accesses the normal vector $n_{TE}$ to a mesh edge $E$ of an element $T \in \Th$ through the keyword \lstinline{normal}.

Parameters can be defined by expressions that contain previously-defined parameters, and functions can access both previously-defined parameters and functions.
Calls between spatial functions use their resolved argument, as in \lstinline{exact_solution_gradient(X)} above (by contrast, inside an integral a bare spatial-function name denotes its value at the current integration point; see, e.g., Listing \ref{lst:load}).
Available geometric contexts in two space dimensions are \lstinline{on element T}, \lstinline{on edge E}, and \lstinline{on edge E of element T}.

Note that the symbols \lstinline{Omega}, \lstinline{T}, \lstinline{F}, \lstinline{E}, and \lstinline{V} are hardwired in the DSL and refer, respectively, to the domain and to a generic element, face, edge, and vertex.

\subsection{Discrete spaces}

For any $Y \in \Mh$, we denote by $\Poly{k}(Y)$ the space spanned by polynomial functions $Y \to \Real$ of total degree up to $k \ge 0$, with the additional convention that $\Poly{-1}(Y) \coloneqq \{0\}$.
The HHO method for problem \eqref{eq:continuous} proposed in \cite{Di-Pietro.Ern.ea:14} is based on the following \emph{discrete space}:
\begin{equation}\label{eq:Uh}
  \underline{U}_h^k \coloneqq \left\{
  \underline{v}_h = ( (v_T)_{T \in \Th}, (v_E)_{E \in \Eh} ) \;:\;
  \text{
    $v_T \in \Poly{k}(T)$ for all $T \in \Th$,
    $v_E \in \Poly{k}(E)$ for all $E \in \Eh$
  }
  \right\}.
\end{equation}
With a small abuse of language, the components of the vector $\underline{v}_h$ are called DOFs in the DSL.
The space \eqref{eq:Uh} is defined by the following block:
\begin{lstlisting}[caption={Discrete space}, label=lst:discrete.space]
space Uh {
  element Poly(k, scalar)
  edge Poly(k, scalar)
}
\end{lstlisting}
Notice that the symbol \lstinline{k}, which is reserved by the DSL and does not need to be declared, refers to the degree of the method as specified by the option \lstinline[style=GenericStyle]{--degree} of the solvers; see Section \ref{sec:solvers}.
\begin{table}\centering
  \caption{Polynomial families that can be used to define DOFs in a discrete space. Here, $Y$ denotes the supporting mesh entity, $T$ an element, and $\mathbb{X}_r$ the scalar, vector, or matrix value space.
    When the rank specification is omitted, it defaults to the one in boldface font.}
  \label{tab:dsl-dof-polynomial-spaces}
  \small
  \renewcommand{\arraystretch}{1.18}
  \setlength{\tabcolsep}{4pt}
  \begin{tabularx}{\textwidth}{@{}
      >{\raggedright\arraybackslash}p{0.30\textwidth}
      >{\raggedright\arraybackslash}p{0.34\textwidth}
      >{\raggedright\arraybackslash}p{0.14\textwidth}
      >{\raggedright\arraybackslash}X@{}}
    \toprule
    \textbf{DSL syntax}
      & \textbf{Polynomial space}
      & \textbf{Rank}
      & \textbf{DOF support} \\
    \midrule

    \lstinline{Poly(m, r)}
      & $\Poly{m}(Y;\mathbb{X}_{r})$
      & \textbf{scalar}, vector, or matrix
      & element, edge, or vertex; on the domain, only
        \lstinline{Poly(0, r)} is allowed \\

    \lstinline{ZeroAveragePoly(m, r)}
      & $\left\{q\in\Poly{m}(Y;\mathbb{X}_{r}) \;:\; \int_Y q=0 \right\}$
      & \textbf{scalar}, vector, or matrix
      & element or edge \\

    \lstinline{GradientPoly(m)} or \lstinline{GradPoly(m)}
      & $\mathcal{G}^{m}(T)
          =\nabla\Poly{m+1}(T)$
      & \textbf{vector}
      & element only \\

    \lstinline{GradientPolyComplement(m)} or \lstinline{GradPolyComplement(m)}
      & $\mathcal{G}^{\mathrm{c},m}(T)
          =(x-x_T)^{\perp}\Poly{m-1}(T)$
      & \textbf{vector}
      & element only \\

    \lstinline{CurlPoly(m)}
      & $\mathcal{R}^{m}(T) = \operatorname{curl} \Poly{m+1}(T)$
      & \textbf{vector}
      & element only \\

    \lstinline{CurlPolyComplement(m, r)}
      & $\mathcal{R}^{\mathrm{c},m}(T) = (x-x_T) \Poly{m-1}(T)$;
     for matrix rank, each row belongs to
     $\mathcal{R}^{\mathrm{c},m}(T)$
      & \textbf{vector} or matrix
      & element only \\


    \lstinline{NedelecPoly(m)}
      & $\mathcal{N}^{m}(T)
          = \big(\mathcal{RT}^{m}(T)\big)^{\perp}
          = \mathcal{G}^{m-1}(T) + \mathcal{G}^{\mathrm{c},m}(T)$
      & \textbf{vector}
      & element only \\

    \lstinline{RaviartThomasPoly(m)}
      & $\mathcal{RT}^{m}(T) = \mathcal{R}^{m-1}(T) + \mathcal{R}^{\mathrm{c},m}(T)$
      & \textbf{vector}
      & element only \\

    \lstinline{orthogonal complement of A relative to B}
      & $B\cap A^{\perp_{L^2(Y)}}$
      & inherited from $A$ and $B$, which must have the same rank
      & element or edge, subject to the restrictions on $A$ and $B$ \\

    \bottomrule
  \end{tabularx}
\end{table}

A list of the main polynomial spaces that can be used to define DOFs in a discrete space is provided in Table \ref{tab:dsl-dof-polynomial-spaces}.
The degree \lstinline{m} can be an integer, \lstinline{k}, \lstinline{k+j}, or \lstinline{k-j}, \lstinline{j} being a constant positive integer.
A negative resolved degree denotes the trivial space.

When implementing more complicated methods, there may be more than one DOF supported by a given mesh entity type.
For this reason, the DSL allows the user to attach a name to a specific DOF, as illustrated in the following variation of the previous code:
\begin{lstlisting}[caption={Discrete space with named DOFs}, label=lst:discrete.space:named.dofs]
space Uh {
  element Poly(k, scalar) called element_dof
  edge Poly(k, scalar) called edge_dof
}
\end{lstlisting}
The DSL also supports the declaration of Cartesian product spaces; see Section \ref{sec:advanced.features:cartesian.product.spaces}.

\subsection{Interpolants}

The meaning of the components of a vector $\underline{v}_h \in \underline{U}_h^k$ is provided by the \emph{interpolation operator} $\underline{I}_h^k : H^1(\Omega) \to \underline{U}_h^k$, which maps a given function $v \in H^1(\Omega)$ onto its \emph{interpolant}
\[
\underline{I}_h^k v \coloneqq ( (\pi_T^k v)_{T \in \Th}, (\pi_E^k v)_{E \in \Eh}) \in \underline{U}_h^k,
\]
where $\pi_Y^k$, $Y \in \Mh$, denotes the $L^2$-orthogonal projector onto $\Poly{k}(Y)$.
When prototyping a polytopal method, it can be useful to interpolate a function, either to compute errors or to test the polynomial consistency of operators, as discussed in the next section.
An interpolant can be defined by an \lstinline{interpolant} block, and a method can contain several interpolants.
The interpolant of the \lstinline{exact_solution} function defined in Listing \ref{lst:parameters.functions} reads
\begin{lstlisting}[caption={Interpolant}, label=lst:interpolant]
interpolant interpolant_exact_solution on Uh {
  on element T:
    dof(T) = l2_project(exact_solution, Poly(k, scalar))
  on edge E:
    dof(E) = l2_project(exact_solution, Poly(k, scalar))
}
\end{lstlisting}
Named DOFs such as those in Listing \ref{lst:discrete.space:named.dofs} can be accessed with the syntax \lstinline{dof(T, element_dof)} and \lstinline{dof(E, edge_dof)}.
Table \ref{tab:interpolant:operations} collects operations that can be used in interpolants relevant to the present work.

\begin{table}
  \caption{Operations available in interpolants.\label{tab:interpolant:operations}}
  \begin{tabularx}{\textwidth}{@{}lX@{}}
    \toprule
    Operation & Required family argument \\
    \midrule
    \term{raviart\_thomas\_interpolate} &
    A \term{RaviartThomasPoly} family of vector rank. \\
    \term{brezzi\_douglas\_marini\_interpolate} &
    A \term{Poly} family of vector rank. \\
    \term{l2\_projection}, \term{l2\_project} &
    Any family accepted at the selected location. \\
    \term{evaluate\_at\_vertex} &
    No family argument. \\
    \bottomrule
  \end{tabularx}
\end{table}

\subsection{Operators}\label{sec:overview:operators}

At the core of polytopal methods are \emph{discrete operators} that reconstruct relevant quantities from the DOFs.
They are perhaps the most original part of the DSL, and what really sets it apart from similar languages developed in the finite element context.

Denote by $\underline{U}_T^k$ the restriction of the space \eqref{eq:Uh} to a mesh element $T \in \Th$, consisting of vectors of the form $\underline{v}_T = (v_T, (v_E)_{E \in \ET})$, where $\ET$ collects the edges of $T$.
We consider here the variation of the HHO method discussed in \cite[Section 4.2]{Di-Pietro.Droniou:20}, where the consistency term is based on the discrete gradient $G_T^k : \underline{U}_T^k \to \Poly{k}(T; \Real^2)$ such that, for all $\underline{v}_T \in \underline{U}_T^k$,
\begin{equation}\label{eq:GT}
  \int_T G_T^k \underline{v}_T \cdot \tau
  = -\int_T v_T \operatorname{div} \tau
  + \sum_{E \in \ET} \int_E v_E (\tau \cdot n_{TE})
  \qquad \forall \tau \in \Poly{k}(T; \Real^2).
\end{equation}
The DSL realization of \eqref{eq:GT} is
\begin{lstlisting}[caption={Gradient reconstruction}, label=lst:gradient.reconstruction]
operator gradient_reconstruction : Uh(v) -> Poly(k, vector) on element T {
  forall tau in Poly(k, vector):
    int(T) gradient_reconstruction(v) dot tau =
    - int(T) dof(v, T) * div(tau)
    + sum_edges(int(E) dof(v, E) * (tau dot normal))
}
\end{lstlisting}
The signature of the operator specifies its domain along with the name of a dummy argument \lstinline{v}, its codomain, and the geometric context (in this case, the element $T$).
As for functions, the geometric contexts \lstinline{on edge E} and \lstinline{on edge E of element T} are also available for two-dimensional problems.
The Galerkin-type variational problem solved by the operator is expressed in natural form, and the language requires that the operator being defined appear in the expression on the left-hand side.
Multiple \lstinline{forall ... in ... :} statements are possible within an operator block; see, e.g., Listing \ref{lst:potential} below.
This is particularly useful when the problem defining the operator is of Petrov--Galerkin type, as for DDR methods.

To build a stabilization in the spirit of \cite[Assumption 2.4]{Di-Pietro.Droniou:20}, we also need the higher-order potential operator $p_T^{k+1} : \underline{U}_T^k \to \Poly{k+1}(T)$ such that, for all $\underline{v}_T \in \underline{U}_T^k$,
\begin{subequations}\label{eq:pT}
  \begin{gather}\label{eq:pT:variational}
    \int_T \nabla p_T^{k+1} \underline{v}_T \cdot \nabla w
    = \int_T G_T^k \underline{v}_T \cdot \nabla w
    \qquad \forall w \in \Poly{k+1}_0(T),
    \\ \label{eq:pT:constraint}
    \int_T p_T^{k+1} \underline{v}_T = \int_T v_T,
  \end{gather}
\end{subequations}
where $\Poly{k+1}_0(T) \coloneqq \left\{ w \in \Poly{k+1}(T) \;:\; \int_T w = 0 \right\}$.
The expression of \eqref{eq:pT} in the language is
\begin{lstlisting}[caption={Potential operator}, label=lst:potential]
operator potential_reconstruction : Uh(u) -> Poly(k+1, scalar) on element T {
  forall w in ZeroAveragePoly(k+1, scalar):
    int(T) grad(potential_reconstruction(u)) dot grad(w) =
    int(T) gradient_reconstruction(u) dot grad(w)
  forall w in Poly(0, scalar):
      int(T) potential_reconstruction(u) * w = int(T) dof(u, T) * w
}
\end{lstlisting}
The DSL contains shortcuts for constraints such as enforcing the average value or a pointwise value for operators defined on edges.
Thus, an equivalent but more compact definition is:
\begin{lstlisting}
operator potential_reconstruction : Uh(u) -> Poly(k+1, scalar) on element T {
  forall w in Poly(k+1, scalar):
    int(T) grad(potential_reconstruction(u)) dot grad(w) =
    int(T) gradient_reconstruction(u) dot grad(w)
  constraint int(T) potential_reconstruction(u) = int(T) dof(u, T)
}
\end{lstlisting}
Notice that, in both cases, we were able to use the previously defined \lstinline{gradient_reconstruction} operator on the right-hand side, closely mimicking \eqref{eq:pT:variational}.

A crucial property of discrete operators is exactness for interpolants of polynomial functions up to a certain degree.
For the operators defined above, we have, in particular,
\begin{equation}\label{eq:polynomial.consistency}
  \text{
    $p_T^{k+1} \underline{I}_T^k v = v$
    and $G_T^k \underline{I}_T^k v = \nabla v$
    for all $v \in \Poly{k+1}(T)$.
  }
\end{equation}
Checking polynomial consistency numerically is a common debugging step.
The DSL in \polyrust offers the possibility to add consistency tests to operators.
Let us assume that linear and quadratic functions of the space coordinates have been defined:
\begin{lstlisting}
function test_linear(vector X) -> scalar = X[0] - 4.0 * X[1] + 7.0
function test_quadratic(vector X) -> scalar = pow(X[0], 2.0) - 3.0 * X[0] * X[1] + pow(X[1], 2.0) - 5.0
\end{lstlisting}
along with the corresponding interpolants \lstinline{interpolant_test_linear} and \lstinline{interpolant_test_quadratic} in a similar way as in Listing \ref{lst:interpolant}.
Consistency tests based on \eqref{eq:polynomial.consistency} can then be added to the \lstinline{potential} operator of Listing \ref{lst:potential} by modifying the corresponding block as follows:
\begin{lstlisting}
operator potential_reconstruction : Uh(u) -> Poly(k+1, scalar) on element T {
  ...
  test exactness for k = 0 against test_linear using interpolant_test_linear
  test exactness for k = 1 against test_quadratic using interpolant_test_quadratic
}
\end{lstlisting}
This syntax may appear elaborate, but it is justified by considering the form of exactness tests for the gradient operator.
Assume that vector-valued functions \lstinline{test_linear_gradient} and \lstinline{test_quadratic_gradient} have been defined to represent the gradients of the above test functions.
An exactness test for the gradient reconstruction would look like this:
\begin{lstlisting}
operator gradient_reconstruction : Uh(v) -> Poly(k, vector) on element T {
  ...
  test exactness for k = 0 against test_linear_gradient using interpolant_test_linear
  test exactness for k = 1 against test_quadratic_gradient using interpolant_test_quadratic
}
\end{lstlisting}
Here, we are telling the solver that it should compare the result of \lstinline{gradient_reconstruction} applied to \lstinline{interpolant_test_linear} (for $k = 0$) or \lstinline{interpolant_test_quadratic} (for $k = 1$) to \lstinline{test_linear_gradient} and \lstinline{test_quadratic_gradient}, respectively.
By running the following command, the details of which will become clear in Section \ref{sec:solvers},
\begin{lstlisting}[style=GenericStyle, caption={Run exactness tests for operators}, label=lst:operator.exactness]
dsl_solver_linear_problem \
  --mesh meshes/unit-square_1.vtk \
  --dsl hho_poisson.dsl \
  --degree 1 \
  --test-operator-exactness 
\end{lstlisting}%
we obtain the following output:
\begin{lstlisting}[style=GenericStyle]
mesh: meshes/unit-square_1.vtk
mesh size = 2.6516504294495546e-1
DSL file = ../polytopal-dsl/dsl/hho_poisson.dsl
operator exactness test degree k = 1
  gradient_reconstruction against test_quadratic_gradient using interpolant_test_quadratic  ok
  potential_reconstruction against test_quadratic using interpolant_test_quadratic          ok
operator exactness summary: 2 passed, 0 failed
\end{lstlisting}
This shows that both consistency tests specified for $k = 1$ passed.

To conclude, we define the high-order difference operators $\delta_T^k \underline{v}_T \coloneqq \pi_T^k(p_T^{k+1} \underline{v}_T - v_T)$ and, for all $E \in \ET$, $\delta_{TE}^k \underline{v}_T \coloneqq \pi_E^k(p_T^{k+1} \underline{v}_T - v_E)$ (see \cite[Eq. (2.19)]{Di-Pietro.Droniou:20}) as follows:
\begin{lstlisting}
operator element_difference : Uh(v) -> Poly(k, scalar) on element T {
  forall q in Poly(k, scalar):
    int(T) element_difference(v) * q =
    int(T) (potential_reconstruction(v) - dof(v, T)) * q
}
  
operator edge_difference : Uh(v) -> Poly(k, scalar) on edge E of element T {
  forall q in Poly(k, scalar):
    int(E) edge_difference(v) * q =
    int(E) (potential_reconstruction(v) - dof(v, E)) * q
}
\end{lstlisting}

\subsection{Forms}

The HHO bilinear form $a_h : \underline{U}_h^k \times \underline{U}_h^k \to \Real$ is defined as follows: for all $(\underline{u}_h, \underline{v}_h) \in \underline{U}_h^k \times \underline{U}_h^k$,
\begin{equation}\label{eq:ah}
  \begin{gathered}
    a_h(\underline{u}_h, \underline{v}_h)
    \coloneqq \sum_{T \in \Th} a_T(\underline{u}_T, \underline{v}_T),
    \\
    a_T(\underline{u}_T, \underline{v}_T)
    \coloneqq \int_T G_T^k \underline{u}_T \cdot G_T^k \underline{v}_T
    + h_T^{-2} \int_T \delta_T^k \underline{u}_T \, \delta_T^k \underline{v}_T
    + h_T^{-1} \sum_{E \in \ET} \int_E \delta_{TE}^k \underline{u}_T \, \delta_{TE}^k \underline{v}_T.
  \end{gathered}
\end{equation}
Here, we have considered the VEM-inspired stabilization discussed in \cite[Example 2.8]{Di-Pietro.Droniou:20} with the element diameter replacing the face diameter everywhere for robustness in the presence of small faces \cite{Droniou.Yemm:22}.
We leave it as an exercise to the reader to figure out the implementation of the classical HHO stabilization of \cite[Example 2.7]{Di-Pietro.Droniou:20}.
The DSL rendition of \eqref{eq:ah} is
\begin{lstlisting}
bilinear form poisson : Uh(trial u) times Uh(test v) {
  sum_elements(
    int(T) gradient_reconstruction(u) dot gradient_reconstruction(v)
    + pow(diameter(T), -2.0) * int(T) element_difference(u) * element_difference(v)
    + pow(diameter(T), -1.0) * int(dT) edge_difference(u) * edge_difference(v)
  )
}
\end{lstlisting}
In the signature, the additional \lstinline{trial} and \lstinline{test} specifications distinguish the trial and test functions, which are associated, respectively, with the columns and rows of the local matrix.
Notice also that we have used the shortcut \lstinline{int(dT) ...} for the commonly used expression \lstinline{sum_edges(int(E) ... )} already encountered in Listing \ref{lst:gradient.reconstruction}.

The DSL provides an analogous syntax for linear forms.
Consider, e.g., the source linear form $\ell_h : \underline{U}_h^k \to \Real$ of the HHO problem such that, for all $\underline{v}_h \in \underline{U}_h^k$,
\[
\ell_h(\underline{v}_h) \coloneqq \sum_{T \in \Th} \int_T f \, v_T.
\]
Its implementation in the DSL reads
\begin{lstlisting}[caption={Load linear form}, label=lst:load]
linear form load : Uh(test v) {
  sum_elements(int(T) source * dof(v, T))
}
\end{lstlisting}

\subsection{Boundary conditions}

Boundary conditions can be enforced strongly by considering the following affine subset of $\underline{U}_h^k$ (cf. \eqref{eq:Uh}) where DOFs associated with boundary edges are appropriately fixed:
\[
\underline{U}_{h,g}^k \coloneqq \left\{
\underline{v}_h \in \underline{U}_h^k \;:\;
\text{$v_E = \pi_E^k g$ for all $E \in \Eh^\dirichlet$}
\right\},
\]
where $\Eh^\dirichlet \coloneqq \left\{ E \in \Eh \;:\; E \subset \Gamma^\dirichlet\right\}$.
Suppose that we want to run a test on a mesh of the unit square where the edges belonging to each side are identified by integer labels, as is usually the case in VTU/VTK files.
We start by defining human-readable labels as follows:
\begin{lstlisting}
boundary labels { left = 1, right = 2, bottom = 3, top = 4 }
\end{lstlisting}
One can then define a \lstinline{boundary conditions} block which, using a syntax similar to that of the interpolant, fixes the value of boundary DOFs:
\begin{lstlisting}
boundary conditions mixed_boundary_conditions on Uh {
  on edge E in left, bottom, top:
    dof(E) = l2_project(exact_solution, Poly(k, scalar))
}
\end{lstlisting}
Notice that here we prescribe the value of boundary DOFs only on $\Gamma^\dirichlet = \text{\lstinline{left}} \cup \text{\lstinline{bottom}} \cup \text{\lstinline{top}}$, while $\Gamma^\neumann = \text{\lstinline{right}}$ contributes a term to the right-hand side of the problem:
\begin{lstlisting}
linear form neumann_contribution : Uh(test v) {
  sum_boundary_edges(right)(int(E) exact_normal_flux * dof(v, E))
}
\end{lstlisting}
For full Dirichlet boundary conditions, it suffices to write \lstinline{on edge E:} instead of \lstinline{on edge E in left, bottom, top:}.

\subsection{Problems}

The HHO problem reads: Find $\underline{u}_h \in \underline{U}_{h,g}^k$ such that
\begin{equation}\label{eq:discrete}
  a_h(\underline{u}_h, \underline{v}_h)
  = \ell_h(\underline{v}_h) - \int_{\Gamma^\neumann} \phi \, v_E
  \qquad \forall \underline{v}_h \in \underline{U}_{h,0}^k,
\end{equation}
where $\underline{U}_{h,0}^k \coloneqq \left\{
\underline{v}_h \in \underline{U}_h^k \;:\;
\text{$v_E = 0$ for all $E \in \Eh^\dirichlet$}
\right\}$.
We now have all the elements to define it in the DSL:
\begin{lstlisting}[caption={Definition of a linear problem}, label=lst:linear.problem]
linear problem hho_poisson_mixed_boundary_conditions_problem on Uh {
  lhs { poisson }
  rhs { load - neumann_contribution }
  boundary conditions mixed_boundary_conditions
}
\end{lstlisting}
Here we define a linear problem whose left-hand side is obtained by assembling only the \lstinline{poisson} bilinear form, while the right-hand side is obtained by assembling the difference between the \lstinline{load} and \lstinline{neumann_contribution} linear forms.
Specifying \lstinline{lhs} and \lstinline{rhs} as algebraic sums of forms helps modularity and favors code reusability.
We could, for instance, envisage an implementation where the consistency and stabilization contributions in the bilinear form \lstinline{poisson} are implemented as separate \lstinline{bilinear form} blocks,
\begin{lstlisting}
bilinear form poisson_consistency : Uh(trial u) times Uh(test v) { ... }
bilinear form poisson_stabilization : Uh(trial u) times Uh(test v) { ... }
\end{lstlisting}
and the \lstinline{lhs} definition in Listing \ref{lst:linear.problem} is modified as follows:
\begin{lstlisting}
lhs { poisson_consistency +  poisson_stabilization }
\end{lstlisting}
In this way, one could, e.g., define multiple stabilization \lstinline{bilinear form} blocks within the same method and select the one to assemble by simply modifying the \lstinline{lhs} line in the \lstinline{linear problem} block.

\subsection{Functionals and error computation}

Once a discrete solution has been obtained, one usually needs to compute quantities of interest or, when running convergence tests to prototype a new numerical method, the error between the discrete solution and its interpolant.
To this end, the DSL offers \emph{functionals}.
Let us consider an example.
Denote by $\underline{u}_h$ the solution of \eqref{eq:discrete} and suppose, for instance, that we want to compute
\begin{equation}\label{eq:error}
  \norm{1,h}{\underline{u}_h - \underline{I}_h^k u},
\end{equation}
where $u$ denotes the weak solution to \eqref{eq:continuous} and the $H^1$-like component norm is such that, for all $\underline{v}_h \in \underline{U}_h^k$,
\begin{equation}\label{eq:norm.1.h}
  \begin{gathered}
    \norm{1,h}{\underline{v}_h} \coloneqq \left(
    \sum_{T \in \Th} \norm{1,T}{\underline{v}_T}^2
    \right)^{\frac12},
    \\
    \norm{1,T}{\underline{v}_T}^2
    \coloneqq \norm{L^2(T)^2}{\nabla v_T}^2
    + h_T^{-1} \sum_{E \in \ET} \norm{L^2(E)}{v_E - v_T}^2.
  \end{gathered}
\end{equation}
The realization of \eqref{eq:norm.1.h} in the DSL reads
\begin{lstlisting}
function h1_component_norm : Uh(v) -> scalar {
  sqrt(
    sum_elements(
                 int(T) squared_norm(grad(dof(v, T)))
                 + pow(diameter(T), -1.0) * int(dT) pow(dof(v, E) - dof(v, T), 2.0)
                )
    )
}
\end{lstlisting}
The definition of functionals is similar to that of functions, except that both the domain (in this case \lstinline{Uh}) and the name of the dummy argument (\lstinline{v} in the example) are specified.
Functionals can also access previously defined forms.
For instance, the functional corresponding to the $H^1$-like norm induced by $a_h$ simply reads
\begin{lstlisting}
function h1_operator_norm : Uh(v) -> scalar {
  sqrt(poisson(v, v))  
}
\end{lstlisting}
We can now request the computation of \eqref{eq:error} in both norms defined above by adding the following lines to the \lstinline{linear problem} block:
\begin{lstlisting}
linear problem hho_poisson_mixed_boundary_conditions_problem on Uh {
  ...
  compute errors using interpolant_exact_solution {
    h1_component_norm, h1_operator_norm
  }
}  
\end{lstlisting}
The \lstinline{report} field followed by a list of functionals instead evaluates them at the discrete solution.

\subsection{Export}

A problem block can be completed with a list of fields to export.
Exportable fields are operators defined on elements.
For the example developed in this section, one could write, for instance,
\begin{lstlisting}
linear problem hho_poisson_mixed_boundary_conditions_problem on Uh {
  ...
  export { potential_reconstruction, potential_field }
}
\end{lstlisting}
where \lstinline{potential_field} is an operator returning the piecewise function $v_h : \Omega \to \Real$ such that $(v_h)_{|T} = v_T$ for all $T \in \Th$:
\begin{lstlisting}
operator potential_field : Uh(u) -> Poly(k, scalar) on element T {
  potential_field(u) = dof(u, T)
}
\end{lstlisting}

\subsection{Numerical examples}

Let us consider two numerical examples that illustrate some of the features discussed above.

\begin{figure}\centering
  \subcaptionbox{HHO}{%
    \begin{minipage}{0.48\textwidth}
\noindent\begin{tikzpicture}
\begin{loglogaxis}[
  width=\linewidth,
  height=0.78\linewidth,
  title={h1 reconstruction norm vs. $h$},
  title style={font=\small},
  label style={font=\scriptsize},
  tick label style={font=\scriptsize},
  grid=both,
  major grid style={gray!40},
  minor grid style={gray!20},
  legend pos=south east,
  legend style={font=\scriptsize, inner xsep=2pt, inner ysep=1pt},
  legend cell align={left}
]
\addplot+[blue, mark=*, thick] coordinates { (5.3033008588991082e-1,1.3852185941045652e0) (2.6516504294495546e-1,7.3480194519079312e-1) (1.3258252147247773e-1,3.7252814774269866e-1) (6.6291260736238949e-2,1.8690747977964822e-1) (3.3145630368119572e-2,9.3534334382056161e-2) };
\addlegendentry{$k = 0$}
\draw[draw=none] (axis cs:5.3033008588991082e-1,1.3852185941045652e0) -- (axis cs:2.6516504294495546e-1,7.3480194519079312e-1) node[midway, above, text=blue, fill=white, fill opacity=0.8, text opacity=1, inner sep=1pt, font=\scriptsize] {0.91};
\draw[draw=none] (axis cs:2.6516504294495546e-1,7.3480194519079312e-1) -- (axis cs:1.3258252147247773e-1,3.7252814774269866e-1) node[midway, above, text=blue, fill=white, fill opacity=0.8, text opacity=1, inner sep=1pt, font=\scriptsize] {0.98};
\draw[draw=none] (axis cs:1.3258252147247773e-1,3.7252814774269866e-1) -- (axis cs:6.6291260736238949e-2,1.8690747977964822e-1) node[midway, above, text=blue, fill=white, fill opacity=0.8, text opacity=1, inner sep=1pt, font=\scriptsize] {1.00};
\draw[draw=none] (axis cs:6.6291260736238949e-2,1.8690747977964822e-1) -- (axis cs:3.3145630368119572e-2,9.3534334382056161e-2) node[midway, above, text=blue, fill=white, fill opacity=0.8, text opacity=1, inner sep=1pt, font=\scriptsize] {1.00};
\addplot+[red, mark=square*, thick] coordinates { (5.3033008588991082e-1,4.2712964834663819e-1) (2.6516504294495546e-1,1.1465038652946545e-1) (1.3258252147247773e-1,2.9631974238373063e-2) (6.6291260736238949e-2,7.5229768510476115e-3) (3.3145630368119572e-2,1.8938563277548170e-3) };
\addlegendentry{$k = 1$}
\draw[draw=none] (axis cs:5.3033008588991082e-1,4.2712964834663819e-1) -- (axis cs:2.6516504294495546e-1,1.1465038652946545e-1) node[midway, above, text=red, fill=white, fill opacity=0.8, text opacity=1, inner sep=1pt, font=\scriptsize] {1.90};
\draw[draw=none] (axis cs:2.6516504294495546e-1,1.1465038652946545e-1) -- (axis cs:1.3258252147247773e-1,2.9631974238373063e-2) node[midway, above, text=red, fill=white, fill opacity=0.8, text opacity=1, inner sep=1pt, font=\scriptsize] {1.95};
\draw[draw=none] (axis cs:1.3258252147247773e-1,2.9631974238373063e-2) -- (axis cs:6.6291260736238949e-2,7.5229768510476115e-3) node[midway, above, text=red, fill=white, fill opacity=0.8, text opacity=1, inner sep=1pt, font=\scriptsize] {1.98};
\draw[draw=none] (axis cs:6.6291260736238949e-2,7.5229768510476115e-3) -- (axis cs:3.3145630368119572e-2,1.8938563277548170e-3) node[midway, above, text=red, fill=white, fill opacity=0.8, text opacity=1, inner sep=1pt, font=\scriptsize] {1.99};
\addplot+[brown, mark=o, thick] coordinates { (5.3033008588991082e-1,6.2489596095899110e-2) (2.6516504294495546e-1,8.9046298124974953e-3) (1.3258252147247773e-1,1.1705956148997691e-3) (6.6291260736238949e-2,1.4833723352783955e-4) (3.3145630368119572e-2,1.8607317975446422e-5) };
\addlegendentry{$k = 2$}
\draw[draw=none] (axis cs:5.3033008588991082e-1,6.2489596095899110e-2) -- (axis cs:2.6516504294495546e-1,8.9046298124974953e-3) node[midway, above, text=brown, fill=white, fill opacity=0.8, text opacity=1, inner sep=1pt, font=\scriptsize] {2.81};
\draw[draw=none] (axis cs:2.6516504294495546e-1,8.9046298124974953e-3) -- (axis cs:1.3258252147247773e-1,1.1705956148997691e-3) node[midway, above, text=brown, fill=white, fill opacity=0.8, text opacity=1, inner sep=1pt, font=\scriptsize] {2.93};
\draw[draw=none] (axis cs:1.3258252147247773e-1,1.1705956148997691e-3) -- (axis cs:6.6291260736238949e-2,1.4833723352783955e-4) node[midway, above, text=brown, fill=white, fill opacity=0.8, text opacity=1, inner sep=1pt, font=\scriptsize] {2.98};
\draw[draw=none] (axis cs:6.6291260736238949e-2,1.4833723352783955e-4) -- (axis cs:3.3145630368119572e-2,1.8607317975446422e-5) node[midway, above, text=brown, fill=white, fill opacity=0.8, text opacity=1, inner sep=1pt, font=\scriptsize] {2.99};
\addplot+[teal, mark=triangle*, thick] coordinates { (5.3033008588991082e-1,1.3536587757983927e-2) (2.6516504294495546e-1,9.2311368524023287e-4) (1.3258252147247773e-1,5.9240605617092894e-5) (6.6291260736238949e-2,3.7474934854642639e-6) (3.3145630368119572e-2,2.3552228272077170e-7) };
\addlegendentry{$k = 3$}
\draw[draw=none] (axis cs:5.3033008588991082e-1,1.3536587757983927e-2) -- (axis cs:2.6516504294495546e-1,9.2311368524023287e-4) node[midway, above, text=teal, fill=white, fill opacity=0.8, text opacity=1, inner sep=1pt, font=\scriptsize] {3.87};
\draw[draw=none] (axis cs:2.6516504294495546e-1,9.2311368524023287e-4) -- (axis cs:1.3258252147247773e-1,5.9240605617092894e-5) node[midway, above, text=teal, fill=white, fill opacity=0.8, text opacity=1, inner sep=1pt, font=\scriptsize] {3.96};
\draw[draw=none] (axis cs:1.3258252147247773e-1,5.9240605617092894e-5) -- (axis cs:6.6291260736238949e-2,3.7474934854642639e-6) node[midway, above, text=teal, fill=white, fill opacity=0.8, text opacity=1, inner sep=1pt, font=\scriptsize] {3.98};
\draw[draw=none] (axis cs:6.6291260736238949e-2,3.7474934854642639e-6) -- (axis cs:3.3145630368119572e-2,2.3552228272077170e-7) node[midway, above, text=teal, fill=white, fill opacity=0.8, text opacity=1, inner sep=1pt, font=\scriptsize] {3.99};
\end{loglogaxis}
\end{tikzpicture}
    \end{minipage}
    \begin{minipage}{0.48\textwidth}
\noindent\begin{tikzpicture}
\begin{loglogaxis}[
  width=\linewidth,
  height=0.78\linewidth,
  title={l2 reconstruction norm vs. $h$},
  title style={font=\small},
  label style={font=\scriptsize},
  tick label style={font=\scriptsize},
  grid=both,
  major grid style={gray!40},
  minor grid style={gray!20},
  legend pos=south east,
  legend style={font=\scriptsize, inner xsep=2pt, inner ysep=1pt},
  legend cell align={left}
]
\addplot+[blue, mark=*, thick] coordinates { (5.3033008588991082e-1,5.6646708908150478e-1) (2.6516504294495546e-1,1.5271238556610581e-1) (1.3258252147247773e-1,3.8864967417023763e-2) (6.6291260736238949e-2,9.7592615633302522e-3) (3.3145630368119572e-2,2.4425066722522700e-3) };
\addlegendentry{$k = 0$}
\draw[draw=none] (axis cs:5.3033008588991082e-1,5.6646708908150478e-1) -- (axis cs:2.6516504294495546e-1,1.5271238556610581e-1) node[midway, above, text=blue, fill=white, fill opacity=0.8, text opacity=1, inner sep=1pt, font=\scriptsize] {1.89};
\draw[draw=none] (axis cs:2.6516504294495546e-1,1.5271238556610581e-1) -- (axis cs:1.3258252147247773e-1,3.8864967417023763e-2) node[midway, above, text=blue, fill=white, fill opacity=0.8, text opacity=1, inner sep=1pt, font=\scriptsize] {1.97};
\draw[draw=none] (axis cs:1.3258252147247773e-1,3.8864967417023763e-2) -- (axis cs:6.6291260736238949e-2,9.7592615633302522e-3) node[midway, above, text=blue, fill=white, fill opacity=0.8, text opacity=1, inner sep=1pt, font=\scriptsize] {1.99};
\draw[draw=none] (axis cs:6.6291260736238949e-2,9.7592615633302522e-3) -- (axis cs:3.3145630368119572e-2,2.4425066722522700e-3) node[midway, above, text=blue, fill=white, fill opacity=0.8, text opacity=1, inner sep=1pt, font=\scriptsize] {2.00};
\addplot+[red, mark=square*, thick] coordinates { (5.3033008588991082e-1,1.9427297180590314e-1) (2.6516504294495546e-1,2.6647821231060669e-2) (1.3258252147247773e-1,3.5189317725379289e-3) (6.6291260736238949e-2,4.5266635271211679e-4) (3.3145630368119572e-2,5.7357632386032103e-5) };
\addlegendentry{$k = 1$}
\draw[draw=none] (axis cs:5.3033008588991082e-1,1.9427297180590314e-1) -- (axis cs:2.6516504294495546e-1,2.6647821231060669e-2) node[midway, above, text=red, fill=white, fill opacity=0.8, text opacity=1, inner sep=1pt, font=\scriptsize] {2.87};
\draw[draw=none] (axis cs:2.6516504294495546e-1,2.6647821231060669e-2) -- (axis cs:1.3258252147247773e-1,3.5189317725379289e-3) node[midway, above, text=red, fill=white, fill opacity=0.8, text opacity=1, inner sep=1pt, font=\scriptsize] {2.92};
\draw[draw=none] (axis cs:1.3258252147247773e-1,3.5189317725379289e-3) -- (axis cs:6.6291260736238949e-2,4.5266635271211679e-4) node[midway, above, text=red, fill=white, fill opacity=0.8, text opacity=1, inner sep=1pt, font=\scriptsize] {2.96};
\draw[draw=none] (axis cs:6.6291260736238949e-2,4.5266635271211679e-4) -- (axis cs:3.3145630368119572e-2,5.7357632386032103e-5) node[midway, above, text=red, fill=white, fill opacity=0.8, text opacity=1, inner sep=1pt, font=\scriptsize] {2.98};
\addplot+[brown, mark=o, thick] coordinates { (5.3033008588991082e-1,2.8536526682477997e-2) (2.6516504294495546e-1,2.0870715417809537e-3) (1.3258252147247773e-1,1.3766473026171800e-4) (6.6291260736238949e-2,8.6993905037404254e-6) (3.3145630368119572e-2,5.4428703764051743e-7) };
\addlegendentry{$k = 2$}
\draw[draw=none] (axis cs:5.3033008588991082e-1,2.8536526682477997e-2) -- (axis cs:2.6516504294495546e-1,2.0870715417809537e-3) node[midway, above, text=brown, fill=white, fill opacity=0.8, text opacity=1, inner sep=1pt, font=\scriptsize] {3.77};
\draw[draw=none] (axis cs:2.6516504294495546e-1,2.0870715417809537e-3) -- (axis cs:1.3258252147247773e-1,1.3766473026171800e-4) node[midway, above, text=brown, fill=white, fill opacity=0.8, text opacity=1, inner sep=1pt, font=\scriptsize] {3.92};
\draw[draw=none] (axis cs:1.3258252147247773e-1,1.3766473026171800e-4) -- (axis cs:6.6291260736238949e-2,8.6993905037404254e-6) node[midway, above, text=brown, fill=white, fill opacity=0.8, text opacity=1, inner sep=1pt, font=\scriptsize] {3.98};
\draw[draw=none] (axis cs:6.6291260736238949e-2,8.6993905037404254e-6) -- (axis cs:3.3145630368119572e-2,5.4428703764051743e-7) node[midway, above, text=brown, fill=white, fill opacity=0.8, text opacity=1, inner sep=1pt, font=\scriptsize] {4.00};
\addplot+[teal, mark=triangle*, thick] coordinates { (5.3033008588991082e-1,6.3165649980844087e-3) (2.6516504294495546e-1,2.2580199098934178e-4) (1.3258252147247773e-1,7.3830401539187628e-6) (6.6291260736238949e-2,2.3531525605692381e-7) (3.3145630368119572e-2,7.4181786922442335e-9) };
\addlegendentry{$k = 3$}
\draw[draw=none] (axis cs:5.3033008588991082e-1,6.3165649980844087e-3) -- (axis cs:2.6516504294495546e-1,2.2580199098934178e-4) node[midway, above, text=teal, fill=white, fill opacity=0.8, text opacity=1, inner sep=1pt, font=\scriptsize] {4.81};
\draw[draw=none] (axis cs:2.6516504294495546e-1,2.2580199098934178e-4) -- (axis cs:1.3258252147247773e-1,7.3830401539187628e-6) node[midway, above, text=teal, fill=white, fill opacity=0.8, text opacity=1, inner sep=1pt, font=\scriptsize] {4.93};
\draw[draw=none] (axis cs:1.3258252147247773e-1,7.3830401539187628e-6) -- (axis cs:6.6291260736238949e-2,2.3531525605692381e-7) node[midway, above, text=teal, fill=white, fill opacity=0.8, text opacity=1, inner sep=1pt, font=\scriptsize] {4.97};
\draw[draw=none] (axis cs:6.6291260736238949e-2,2.3531525605692381e-7) -- (axis cs:3.3145630368119572e-2,7.4181786922442335e-9) node[midway, above, text=teal, fill=white, fill opacity=0.8, text opacity=1, inner sep=1pt, font=\scriptsize] {4.99};
\end{loglogaxis}
\end{tikzpicture}
    \end{minipage}
  }
  \\
  \subcaptionbox{DDR}{%
    \begin{minipage}{0.48\textwidth}
\noindent\begin{tikzpicture}
\begin{loglogaxis}[
  width=\linewidth,
  height=0.78\linewidth,
  title={h1 potential norm vs. $h$},
  title style={font=\small},
  label style={font=\scriptsize},
  tick label style={font=\scriptsize},
  grid=both,
  major grid style={gray!40},
  minor grid style={gray!20},
  legend pos=south east,
  legend style={font=\scriptsize, inner xsep=2pt, inner ysep=1pt},
  legend cell align={left}
]
\addplot+[blue, mark=*, thick] coordinates { (5.3033008588991082e-1,9.3691856643813465e-2) (2.6516504294495546e-1,3.8509630509489017e-2) (1.3258252147247773e-1,1.2526681362319238e-2) (6.6291260736238949e-2,3.6668930298708259e-3) (3.3145630368119572e-2,1.0246905333972322e-3) };
\addlegendentry{$k = 0$}
\draw[draw=none] (axis cs:5.3033008588991082e-1,9.3691856643813465e-2) -- (axis cs:2.6516504294495546e-1,3.8509630509489017e-2) node[midway, above, text=blue, fill=white, fill opacity=0.8, text opacity=1, inner sep=1pt, font=\scriptsize] {1.28};
\draw[draw=none] (axis cs:2.6516504294495546e-1,3.8509630509489017e-2) -- (axis cs:1.3258252147247773e-1,1.2526681362319238e-2) node[midway, above, text=blue, fill=white, fill opacity=0.8, text opacity=1, inner sep=1pt, font=\scriptsize] {1.62};
\draw[draw=none] (axis cs:1.3258252147247773e-1,1.2526681362319238e-2) -- (axis cs:6.6291260736238949e-2,3.6668930298708259e-3) node[midway, above, text=blue, fill=white, fill opacity=0.8, text opacity=1, inner sep=1pt, font=\scriptsize] {1.77};
\draw[draw=none] (axis cs:6.6291260736238949e-2,3.6668930298708259e-3) -- (axis cs:3.3145630368119572e-2,1.0246905333972322e-3) node[midway, above, text=blue, fill=white, fill opacity=0.8, text opacity=1, inner sep=1pt, font=\scriptsize] {1.84};
\addplot+[red, mark=square*, thick] coordinates { (5.3033008588991082e-1,1.3382686046624634e-1) (2.6516504294495546e-1,2.7441544456902913e-2) (1.3258252147247773e-1,4.8728272309979531e-3) (6.6291260736238949e-2,8.6511763320104009e-4) (3.3145630368119572e-2,1.6634295079400528e-4) };
\addlegendentry{$k = 1$}
\draw[draw=none] (axis cs:5.3033008588991082e-1,1.3382686046624634e-1) -- (axis cs:2.6516504294495546e-1,2.7441544456902913e-2) node[midway, above, text=red, fill=white, fill opacity=0.8, text opacity=1, inner sep=1pt, font=\scriptsize] {2.29};
\draw[draw=none] (axis cs:2.6516504294495546e-1,2.7441544456902913e-2) -- (axis cs:1.3258252147247773e-1,4.8728272309979531e-3) node[midway, above, text=red, fill=white, fill opacity=0.8, text opacity=1, inner sep=1pt, font=\scriptsize] {2.49};
\draw[draw=none] (axis cs:1.3258252147247773e-1,4.8728272309979531e-3) -- (axis cs:6.6291260736238949e-2,8.6511763320104009e-4) node[midway, above, text=red, fill=white, fill opacity=0.8, text opacity=1, inner sep=1pt, font=\scriptsize] {2.49};
\draw[draw=none] (axis cs:6.6291260736238949e-2,8.6511763320104009e-4) -- (axis cs:3.3145630368119572e-2,1.6634295079400528e-4) node[midway, above, text=red, fill=white, fill opacity=0.8, text opacity=1, inner sep=1pt, font=\scriptsize] {2.38};
\addplot+[brown, mark=o, thick] coordinates { (5.3033008588991082e-1,6.9785458210006200e-3) (2.6516504294495546e-1,1.1494823596417730e-3) (1.3258252147247773e-1,1.5020425844918675e-4) (6.6291260736238949e-2,1.9059055593675081e-5) (3.3145630368119572e-2,2.3954417873509520e-6) };
\addlegendentry{$k = 2$}
\draw[draw=none] (axis cs:5.3033008588991082e-1,6.9785458210006200e-3) -- (axis cs:2.6516504294495546e-1,1.1494823596417730e-3) node[midway, above, text=brown, fill=white, fill opacity=0.8, text opacity=1, inner sep=1pt, font=\scriptsize] {2.60};
\draw[draw=none] (axis cs:2.6516504294495546e-1,1.1494823596417730e-3) -- (axis cs:1.3258252147247773e-1,1.5020425844918675e-4) node[midway, above, text=brown, fill=white, fill opacity=0.8, text opacity=1, inner sep=1pt, font=\scriptsize] {2.94};
\draw[draw=none] (axis cs:1.3258252147247773e-1,1.5020425844918675e-4) -- (axis cs:6.6291260736238949e-2,1.9059055593675081e-5) node[midway, above, text=brown, fill=white, fill opacity=0.8, text opacity=1, inner sep=1pt, font=\scriptsize] {2.98};
\draw[draw=none] (axis cs:6.6291260736238949e-2,1.9059055593675081e-5) -- (axis cs:3.3145630368119572e-2,2.3954417873509520e-6) node[midway, above, text=brown, fill=white, fill opacity=0.8, text opacity=1, inner sep=1pt, font=\scriptsize] {2.99};
\addplot+[teal, mark=triangle*, thick] coordinates { (5.3033008588991082e-1,2.3299882261433782e-3) (2.6516504294495546e-1,1.3319520846108879e-4) (1.3258252147247773e-1,7.2924843202491089e-6) (6.6291260736238949e-2,4.1495526554259033e-7) (3.3145630368119572e-2,2.4617722375702347e-8) };
\addlegendentry{$k = 3$}
\draw[draw=none] (axis cs:5.3033008588991082e-1,2.3299882261433782e-3) -- (axis cs:2.6516504294495546e-1,1.3319520846108879e-4) node[midway, above, text=teal, fill=white, fill opacity=0.8, text opacity=1, inner sep=1pt, font=\scriptsize] {4.13};
\draw[draw=none] (axis cs:2.6516504294495546e-1,1.3319520846108879e-4) -- (axis cs:1.3258252147247773e-1,7.2924843202491089e-6) node[midway, above, text=teal, fill=white, fill opacity=0.8, text opacity=1, inner sep=1pt, font=\scriptsize] {4.19};
\draw[draw=none] (axis cs:1.3258252147247773e-1,7.2924843202491089e-6) -- (axis cs:6.6291260736238949e-2,4.1495526554259033e-7) node[midway, above, text=teal, fill=white, fill opacity=0.8, text opacity=1, inner sep=1pt, font=\scriptsize] {4.14};
\draw[draw=none] (axis cs:6.6291260736238949e-2,4.1495526554259033e-7) -- (axis cs:3.3145630368119572e-2,2.4617722375702347e-8) node[midway, above, text=teal, fill=white, fill opacity=0.8, text opacity=1, inner sep=1pt, font=\scriptsize] {4.08};
\end{loglogaxis}
\end{tikzpicture}
    \end{minipage}
    \begin{minipage}{0.48\textwidth}
\noindent\begin{tikzpicture}
\begin{loglogaxis}[
  width=\linewidth,
  height=0.78\linewidth,
  title={l2 potential norm vs. $h$},
  title style={font=\small},
  label style={font=\scriptsize},
  tick label style={font=\scriptsize},
  grid=both,
  major grid style={gray!40},
  minor grid style={gray!20},
  legend pos=south east,
  legend style={font=\scriptsize, inner xsep=2pt, inner ysep=1pt},
  legend cell align={left}
]
\addplot+[blue, mark=*, thick] coordinates { (5.3033008588991082e-1,1.7998241351909298e-2) (2.6516504294495546e-1,4.3337629720627763e-3) (1.3258252147247773e-1,1.1224917332642425e-3) (6.6291260736238949e-2,2.8652682934108214e-4) (3.3145630368119572e-2,7.2231132105402701e-5) };
\addlegendentry{$k = 0$}
\draw[draw=none] (axis cs:5.3033008588991082e-1,1.7998241351909298e-2) -- (axis cs:2.6516504294495546e-1,4.3337629720627763e-3) node[midway, above, text=blue, fill=white, fill opacity=0.8, text opacity=1, inner sep=1pt, font=\scriptsize] {2.05};
\draw[draw=none] (axis cs:2.6516504294495546e-1,4.3337629720627763e-3) -- (axis cs:1.3258252147247773e-1,1.1224917332642425e-3) node[midway, above, text=blue, fill=white, fill opacity=0.8, text opacity=1, inner sep=1pt, font=\scriptsize] {1.95};
\draw[draw=none] (axis cs:1.3258252147247773e-1,1.1224917332642425e-3) -- (axis cs:6.6291260736238949e-2,2.8652682934108214e-4) node[midway, above, text=blue, fill=white, fill opacity=0.8, text opacity=1, inner sep=1pt, font=\scriptsize] {1.97};
\draw[draw=none] (axis cs:6.6291260736238949e-2,2.8652682934108214e-4) -- (axis cs:3.3145630368119572e-2,7.2231132105402701e-5) node[midway, above, text=blue, fill=white, fill opacity=0.8, text opacity=1, inner sep=1pt, font=\scriptsize] {1.99};
\addplot+[red, mark=square*, thick] coordinates { (5.3033008588991082e-1,2.6396787187047532e-2) (2.6516504294495546e-1,3.8092504209129257e-3) (1.3258252147247773e-1,3.8613975796522950e-4) (6.6291260736238949e-2,3.3214084856489041e-5) (3.3145630368119572e-2,2.8164868213013659e-6) };
\addlegendentry{$k = 1$}
\draw[draw=none] (axis cs:5.3033008588991082e-1,2.6396787187047532e-2) -- (axis cs:2.6516504294495546e-1,3.8092504209129257e-3) node[midway, above, text=red, fill=white, fill opacity=0.8, text opacity=1, inner sep=1pt, font=\scriptsize] {2.79};
\draw[draw=none] (axis cs:2.6516504294495546e-1,3.8092504209129257e-3) -- (axis cs:1.3258252147247773e-1,3.8613975796522950e-4) node[midway, above, text=red, fill=white, fill opacity=0.8, text opacity=1, inner sep=1pt, font=\scriptsize] {3.30};
\draw[draw=none] (axis cs:1.3258252147247773e-1,3.8613975796522950e-4) -- (axis cs:6.6291260736238949e-2,3.3214084856489041e-5) node[midway, above, text=red, fill=white, fill opacity=0.8, text opacity=1, inner sep=1pt, font=\scriptsize] {3.54};
\draw[draw=none] (axis cs:6.6291260736238949e-2,3.3214084856489041e-5) -- (axis cs:3.3145630368119572e-2,2.8164868213013659e-6) node[midway, above, text=red, fill=white, fill opacity=0.8, text opacity=1, inner sep=1pt, font=\scriptsize] {3.56};
\addplot+[brown, mark=o, thick] coordinates { (5.3033008588991082e-1,5.8058582458241732e-4) (2.6516504294495546e-1,5.3496220878783251e-5) (1.3258252147247773e-1,2.8122030200788654e-6) (6.6291260736238949e-2,1.5653456950769756e-7) (3.3145630368119572e-2,9.1618205076024965e-9) };
\addlegendentry{$k = 2$}
\draw[draw=none] (axis cs:5.3033008588991082e-1,5.8058582458241732e-4) -- (axis cs:2.6516504294495546e-1,5.3496220878783251e-5) node[midway, above, text=brown, fill=white, fill opacity=0.8, text opacity=1, inner sep=1pt, font=\scriptsize] {3.44};
\draw[draw=none] (axis cs:2.6516504294495546e-1,5.3496220878783251e-5) -- (axis cs:1.3258252147247773e-1,2.8122030200788654e-6) node[midway, above, text=brown, fill=white, fill opacity=0.8, text opacity=1, inner sep=1pt, font=\scriptsize] {4.25};
\draw[draw=none] (axis cs:1.3258252147247773e-1,2.8122030200788654e-6) -- (axis cs:6.6291260736238949e-2,1.5653456950769756e-7) node[midway, above, text=brown, fill=white, fill opacity=0.8, text opacity=1, inner sep=1pt, font=\scriptsize] {4.17};
\draw[draw=none] (axis cs:6.6291260736238949e-2,1.5653456950769756e-7) -- (axis cs:3.3145630368119572e-2,9.1618205076024965e-9) node[midway, above, text=brown, fill=white, fill opacity=0.8, text opacity=1, inner sep=1pt, font=\scriptsize] {4.09};
\addplot+[teal, mark=triangle*, thick] coordinates { (5.3033008588991082e-1,4.0199878777502781e-4) (2.6516504294495546e-1,1.3967754498004528e-5) (1.3258252147247773e-1,3.6489069244419239e-7) (6.6291260736238949e-2,8.3841520050612316e-9) (3.3145630368119572e-2,1.9073196635729773e-10) };
\addlegendentry{$k = 3$}
\draw[draw=none] (axis cs:5.3033008588991082e-1,4.0199878777502781e-4) -- (axis cs:2.6516504294495546e-1,1.3967754498004528e-5) node[midway, above, text=teal, fill=white, fill opacity=0.8, text opacity=1, inner sep=1pt, font=\scriptsize] {4.85};
\draw[draw=none] (axis cs:2.6516504294495546e-1,1.3967754498004528e-5) -- (axis cs:1.3258252147247773e-1,3.6489069244419239e-7) node[midway, above, text=teal, fill=white, fill opacity=0.8, text opacity=1, inner sep=1pt, font=\scriptsize] {5.26};
\draw[draw=none] (axis cs:1.3258252147247773e-1,3.6489069244419239e-7) -- (axis cs:6.6291260736238949e-2,8.3841520050612316e-9) node[midway, above, text=teal, fill=white, fill opacity=0.8, text opacity=1, inner sep=1pt, font=\scriptsize] {5.44};
\draw[draw=none] (axis cs:6.6291260736238949e-2,8.3841520050612316e-9) -- (axis cs:3.3145630368119572e-2,1.9073196635729773e-10) node[midway, above, text=teal, fill=white, fill opacity=0.8, text opacity=1, inner sep=1pt, font=\scriptsize] {5.46};
\end{loglogaxis}
\end{tikzpicture}
    \end{minipage}
  }
  \caption{Comparison of HHO (above) and nodal DDR (below) solutions of the Poisson problem on a refined sequence of triangular meshes.}
  \label{fig:hho-ddr-poisson-comparison:tria}
\end{figure}
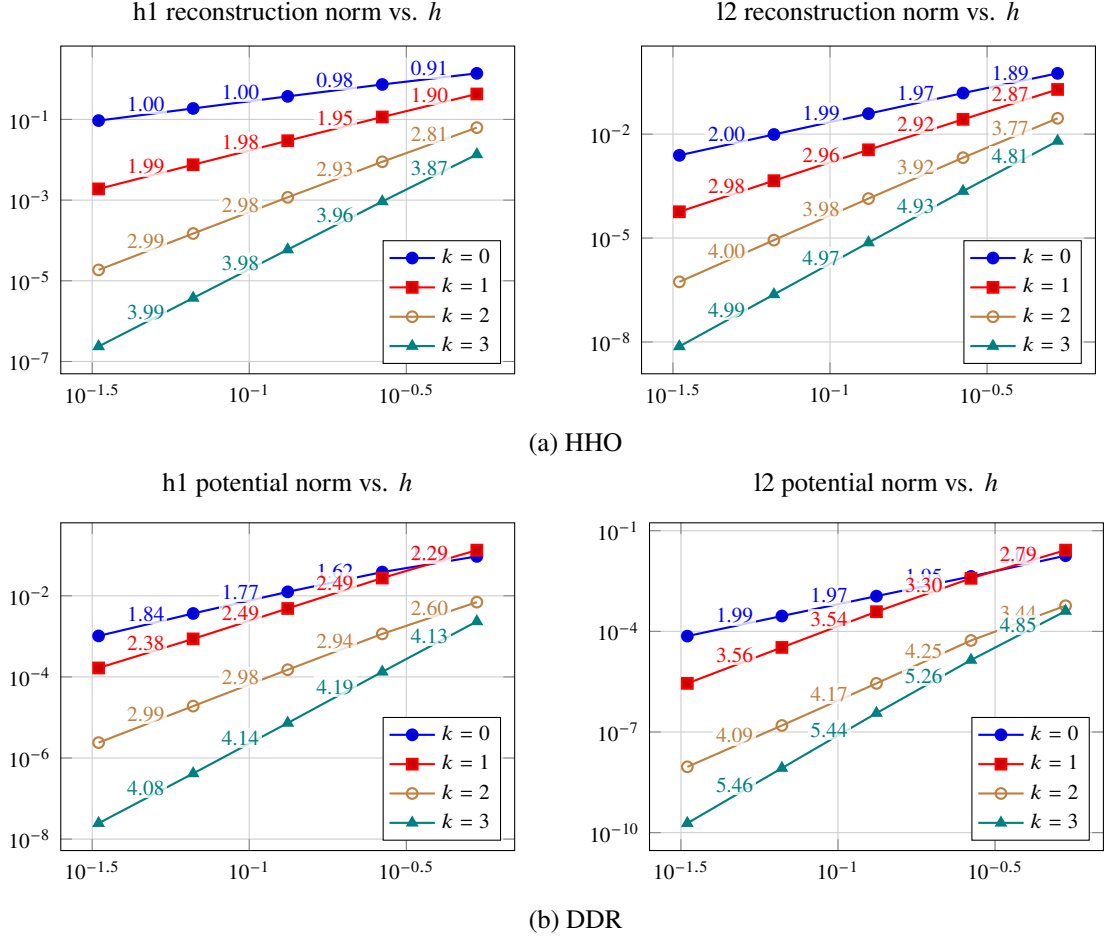

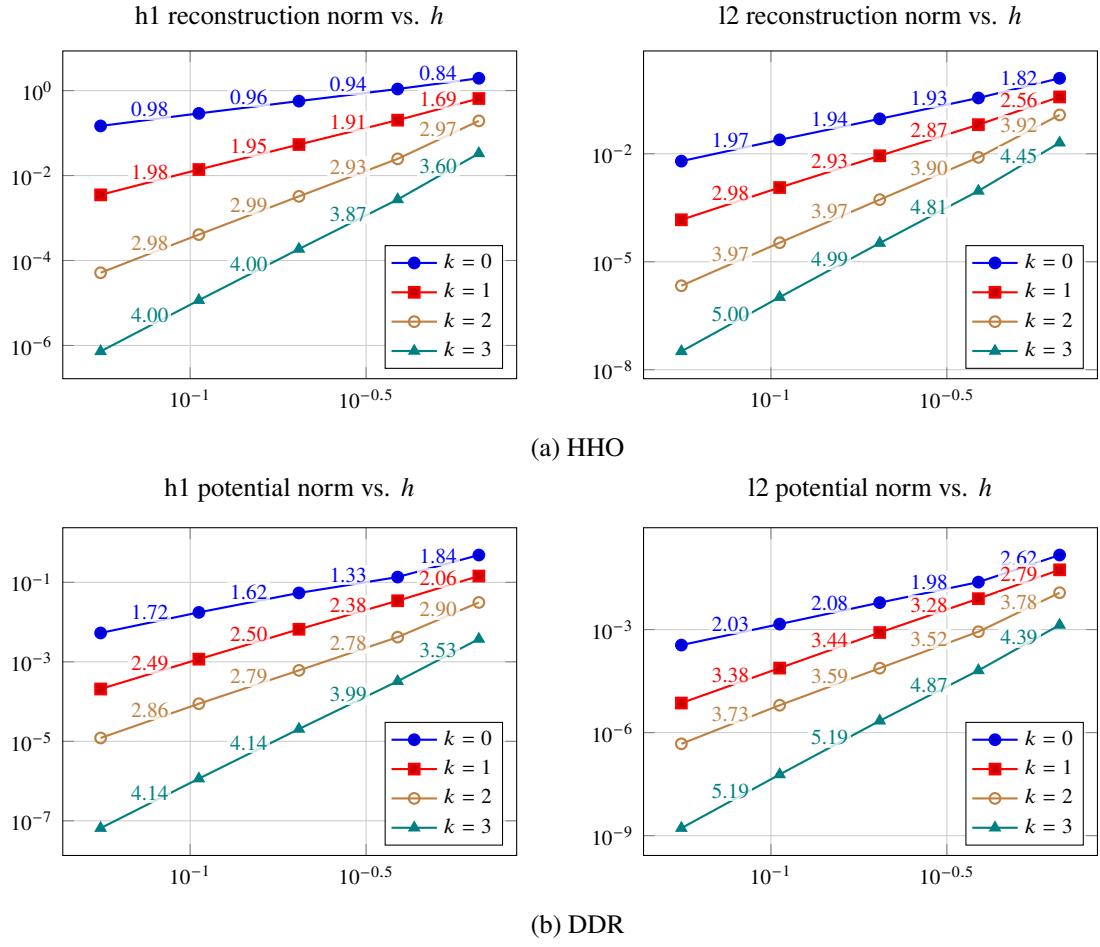
\begin{figure}\centering
  \subcaptionbox{HHO}{%
    \begin{minipage}{0.48\textwidth}
\noindent\begin{tikzpicture}
\begin{loglogaxis}[
  width=\linewidth,
  height=0.78\linewidth,
  title={h1 reconstruction norm vs. $h$},
  title style={font=\small},
  label style={font=\scriptsize},
  tick label style={font=\scriptsize},
  grid=both,
  major grid style={gray!40},
  minor grid style={gray!20},
  legend pos=south east,
  legend style={font=\scriptsize, inner xsep=2pt, inner ysep=1pt},
  legend cell align={left}
]
\addplot+[blue, mark=*, thick] coordinates { (6.6146214270788017e-1,1.9520579261361635e0) (3.8885226009961682e-1,1.0898881703897212e0) (2.0379771530539792e-1,5.6616641708103743e-1) (1.0586955040260637e-1,2.9067224144265275e-1) (5.5649832160483637e-2,1.4719623064468335e-1) };
\addlegendentry{$k = 0$}
\draw[draw=none] (axis cs:6.6146214270788017e-1,1.9520579261361635e0) -- (axis cs:3.8885226009961682e-1,1.0898881703897212e0) node[midway, above, text=blue, fill=white, fill opacity=0.8, text opacity=1, inner sep=1pt, font=\scriptsize] {0.84};
\draw[draw=none] (axis cs:3.8885226009961682e-1,1.0898881703897212e0) -- (axis cs:2.0379771530539792e-1,5.6616641708103743e-1) node[midway, above, text=blue, fill=white, fill opacity=0.8, text opacity=1, inner sep=1pt, font=\scriptsize] {0.94};
\draw[draw=none] (axis cs:2.0379771530539792e-1,5.6616641708103743e-1) -- (axis cs:1.0586955040260637e-1,2.9067224144265275e-1) node[midway, above, text=blue, fill=white, fill opacity=0.8, text opacity=1, inner sep=1pt, font=\scriptsize] {0.96};
\draw[draw=none] (axis cs:1.0586955040260637e-1,2.9067224144265275e-1) -- (axis cs:5.5649832160483637e-2,1.4719623064468335e-1) node[midway, above, text=blue, fill=white, fill opacity=0.8, text opacity=1, inner sep=1pt, font=\scriptsize] {0.98};
\addplot+[red, mark=square*, thick] coordinates { (6.6146214270788017e-1,6.5135900351137210e-1) (3.8885226009961682e-1,2.0147254402243656e-1) (2.0379771530539792e-1,5.3743421240047121e-2) (1.0586955040260637e-1,1.3881769668608076e-2) (5.5649832160483637e-2,3.5134886700665564e-3) };
\addlegendentry{$k = 1$}
\draw[draw=none] (axis cs:6.6146214270788017e-1,6.5135900351137210e-1) -- (axis cs:3.8885226009961682e-1,2.0147254402243656e-1) node[midway, above, text=red, fill=white, fill opacity=0.8, text opacity=1, inner sep=1pt, font=\scriptsize] {1.69};
\draw[draw=none] (axis cs:3.8885226009961682e-1,2.0147254402243656e-1) -- (axis cs:2.0379771530539792e-1,5.3743421240047121e-2) node[midway, above, text=red, fill=white, fill opacity=0.8, text opacity=1, inner sep=1pt, font=\scriptsize] {1.91};
\draw[draw=none] (axis cs:2.0379771530539792e-1,5.3743421240047121e-2) -- (axis cs:1.0586955040260637e-1,1.3881769668608076e-2) node[midway, above, text=red, fill=white, fill opacity=0.8, text opacity=1, inner sep=1pt, font=\scriptsize] {1.95};
\draw[draw=none] (axis cs:1.0586955040260637e-1,1.3881769668608076e-2) -- (axis cs:5.5649832160483637e-2,3.5134886700665564e-3) node[midway, above, text=red, fill=white, fill opacity=0.8, text opacity=1, inner sep=1pt, font=\scriptsize] {1.98};
\addplot+[brown, mark=o, thick] coordinates { (6.6146214270788017e-1,1.9374012544073266e-1) (3.8885226009961682e-1,2.4735941990417166e-2) (2.0379771530539792e-1,3.2365315553319935e-3) (1.0586955040260637e-1,4.0739107097389338e-4) (5.5649832160483637e-2,5.1498633987435318e-5) };
\addlegendentry{$k = 2$}
\draw[draw=none] (axis cs:6.6146214270788017e-1,1.9374012544073266e-1) -- (axis cs:3.8885226009961682e-1,2.4735941990417166e-2) node[midway, above, text=brown, fill=white, fill opacity=0.8, text opacity=1, inner sep=1pt, font=\scriptsize] {2.97};
\draw[draw=none] (axis cs:3.8885226009961682e-1,2.4735941990417166e-2) -- (axis cs:2.0379771530539792e-1,3.2365315553319935e-3) node[midway, above, text=brown, fill=white, fill opacity=0.8, text opacity=1, inner sep=1pt, font=\scriptsize] {2.93};
\draw[draw=none] (axis cs:2.0379771530539792e-1,3.2365315553319935e-3) -- (axis cs:1.0586955040260637e-1,4.0739107097389338e-4) node[midway, above, text=brown, fill=white, fill opacity=0.8, text opacity=1, inner sep=1pt, font=\scriptsize] {2.99};
\draw[draw=none] (axis cs:1.0586955040260637e-1,4.0739107097389338e-4) -- (axis cs:5.5649832160483637e-2,5.1498633987435318e-5) node[midway, above, text=brown, fill=white, fill opacity=0.8, text opacity=1, inner sep=1pt, font=\scriptsize] {2.98};
\addplot+[teal, mark=triangle*, thick] coordinates { (6.6146214270788017e-1,3.3188873597905968e-2) (3.8885226009961682e-1,2.7330909862667324e-3) (2.0379771530539792e-1,1.8639279650443116e-4) (1.0586955040260637e-1,1.1617125683980934e-5) (5.5649832160483637e-2,7.2511357325707996e-7) };
\addlegendentry{$k = 3$}
\draw[draw=none] (axis cs:6.6146214270788017e-1,3.3188873597905968e-2) -- (axis cs:3.8885226009961682e-1,2.7330909862667324e-3) node[midway, above, text=teal, fill=white, fill opacity=0.8, text opacity=1, inner sep=1pt, font=\scriptsize] {3.60};
\draw[draw=none] (axis cs:3.8885226009961682e-1,2.7330909862667324e-3) -- (axis cs:2.0379771530539792e-1,1.8639279650443116e-4) node[midway, above, text=teal, fill=white, fill opacity=0.8, text opacity=1, inner sep=1pt, font=\scriptsize] {3.87};
\draw[draw=none] (axis cs:2.0379771530539792e-1,1.8639279650443116e-4) -- (axis cs:1.0586955040260637e-1,1.1617125683980934e-5) node[midway, above, text=teal, fill=white, fill opacity=0.8, text opacity=1, inner sep=1pt, font=\scriptsize] {4.00};
\draw[draw=none] (axis cs:1.0586955040260637e-1,1.1617125683980934e-5) -- (axis cs:5.5649832160483637e-2,7.2511357325707996e-7) node[midway, above, text=teal, fill=white, fill opacity=0.8, text opacity=1, inner sep=1pt, font=\scriptsize] {4.00};
\end{loglogaxis}
\end{tikzpicture}
    \end{minipage}
    \begin{minipage}{0.48\textwidth}
\noindent\begin{tikzpicture}
\begin{loglogaxis}[
  width=\linewidth,
  height=0.78\linewidth,
  title={l2 reconstruction norm vs. $h$},
  title style={font=\small},
  label style={font=\scriptsize},
  tick label style={font=\scriptsize},
  grid=both,
  major grid style={gray!40},
  minor grid style={gray!20},
  legend pos=south east,
  legend style={font=\scriptsize, inner xsep=2pt, inner ysep=1pt},
  legend cell align={left}
]
\addplot+[blue, mark=*, thick] coordinates { (6.6146214270788017e-1,1.2538295714196024e0) (3.8885226009961682e-1,3.5574981188701338e-1) (2.0379771530539792e-1,9.3542708671342262e-2) (1.0586955040260637e-1,2.4432335793777963e-2) (5.5649832160483637e-2,6.2548277924647190e-3) };
\addlegendentry{$k = 0$}
\draw[draw=none] (axis cs:6.6146214270788017e-1,1.2538295714196024e0) -- (axis cs:3.8885226009961682e-1,3.5574981188701338e-1) node[midway, above, text=blue, fill=white, fill opacity=0.8, text opacity=1, inner sep=1pt, font=\scriptsize] {1.82};
\draw[draw=none] (axis cs:3.8885226009961682e-1,3.5574981188701338e-1) -- (axis cs:2.0379771530539792e-1,9.3542708671342262e-2) node[midway, above, text=blue, fill=white, fill opacity=0.8, text opacity=1, inner sep=1pt, font=\scriptsize] {1.93};
\draw[draw=none] (axis cs:2.0379771530539792e-1,9.3542708671342262e-2) -- (axis cs:1.0586955040260637e-1,2.4432335793777963e-2) node[midway, above, text=blue, fill=white, fill opacity=0.8, text opacity=1, inner sep=1pt, font=\scriptsize] {1.94};
\draw[draw=none] (axis cs:1.0586955040260637e-1,2.4432335793777963e-2) -- (axis cs:5.5649832160483637e-2,6.2548277924647190e-3) node[midway, above, text=blue, fill=white, fill opacity=0.8, text opacity=1, inner sep=1pt, font=\scriptsize] {1.97};
\addplot+[red, mark=square*, thick] coordinates { (6.6146214270788017e-1,3.8030953910935600e-1) (3.8885226009961682e-1,6.4465986689981367e-2) (2.0379771530539792e-1,8.8081586345423145e-3) (1.0586955040260637e-1,1.1541600892693046e-3) (5.5649832160483637e-2,1.4676461787680389e-4) };
\addlegendentry{$k = 1$}
\draw[draw=none] (axis cs:6.6146214270788017e-1,3.8030953910935600e-1) -- (axis cs:3.8885226009961682e-1,6.4465986689981367e-2) node[midway, above, text=red, fill=white, fill opacity=0.8, text opacity=1, inner sep=1pt, font=\scriptsize] {2.56};
\draw[draw=none] (axis cs:3.8885226009961682e-1,6.4465986689981367e-2) -- (axis cs:2.0379771530539792e-1,8.8081586345423145e-3) node[midway, above, text=red, fill=white, fill opacity=0.8, text opacity=1, inner sep=1pt, font=\scriptsize] {2.87};
\draw[draw=none] (axis cs:2.0379771530539792e-1,8.8081586345423145e-3) -- (axis cs:1.0586955040260637e-1,1.1541600892693046e-3) node[midway, above, text=red, fill=white, fill opacity=0.8, text opacity=1, inner sep=1pt, font=\scriptsize] {2.93};
\draw[draw=none] (axis cs:1.0586955040260637e-1,1.1541600892693046e-3) -- (axis cs:5.5649832160483637e-2,1.4676461787680389e-4) node[midway, above, text=red, fill=white, fill opacity=0.8, text opacity=1, inner sep=1pt, font=\scriptsize] {2.98};
\addplot+[brown, mark=o, thick] coordinates { (6.6146214270788017e-1,1.2033971933932998e-1) (3.8885226009961682e-1,7.9526527268324824e-3) (2.0379771530539792e-1,5.3189187731609449e-4) (1.0586955040260637e-1,3.3845534246453790e-5) (5.5649832160483637e-2,2.1571671949367520e-6) };
\addlegendentry{$k = 2$}
\draw[draw=none] (axis cs:6.6146214270788017e-1,1.2033971933932998e-1) -- (axis cs:3.8885226009961682e-1,7.9526527268324824e-3) node[midway, above, text=brown, fill=white, fill opacity=0.8, text opacity=1, inner sep=1pt, font=\scriptsize] {3.92};
\draw[draw=none] (axis cs:3.8885226009961682e-1,7.9526527268324824e-3) -- (axis cs:2.0379771530539792e-1,5.3189187731609449e-4) node[midway, above, text=brown, fill=white, fill opacity=0.8, text opacity=1, inner sep=1pt, font=\scriptsize] {3.90};
\draw[draw=none] (axis cs:2.0379771530539792e-1,5.3189187731609449e-4) -- (axis cs:1.0586955040260637e-1,3.3845534246453790e-5) node[midway, above, text=brown, fill=white, fill opacity=0.8, text opacity=1, inner sep=1pt, font=\scriptsize] {3.97};
\draw[draw=none] (axis cs:1.0586955040260637e-1,3.3845534246453790e-5) -- (axis cs:5.5649832160483637e-2,2.1571671949367520e-6) node[midway, above, text=brown, fill=white, fill opacity=0.8, text opacity=1, inner sep=1pt, font=\scriptsize] {3.97};
\addplot+[teal, mark=triangle*, thick] coordinates { (6.6146214270788017e-1,2.0238148298103050e-2) (3.8885226009961682e-1,9.2803075631356813e-4) (2.0379771530539792e-1,3.2999763241671206e-5) (1.0586955040260637e-1,1.0362609280036809e-6) (5.5649832160483637e-2,3.2352244669103991e-8) };
\addlegendentry{$k = 3$}
\draw[draw=none] (axis cs:6.6146214270788017e-1,2.0238148298103050e-2) -- (axis cs:3.8885226009961682e-1,9.2803075631356813e-4) node[midway, above, text=teal, fill=white, fill opacity=0.8, text opacity=1, inner sep=1pt, font=\scriptsize] {4.45};
\draw[draw=none] (axis cs:3.8885226009961682e-1,9.2803075631356813e-4) -- (axis cs:2.0379771530539792e-1,3.2999763241671206e-5) node[midway, above, text=teal, fill=white, fill opacity=0.8, text opacity=1, inner sep=1pt, font=\scriptsize] {4.81};
\draw[draw=none] (axis cs:2.0379771530539792e-1,3.2999763241671206e-5) -- (axis cs:1.0586955040260637e-1,1.0362609280036809e-6) node[midway, above, text=teal, fill=white, fill opacity=0.8, text opacity=1, inner sep=1pt, font=\scriptsize] {4.99};
\draw[draw=none] (axis cs:1.0586955040260637e-1,1.0362609280036809e-6) -- (axis cs:5.5649832160483637e-2,3.2352244669103991e-8) node[midway, above, text=teal, fill=white, fill opacity=0.8, text opacity=1, inner sep=1pt, font=\scriptsize] {5.00};
\end{loglogaxis}
\end{tikzpicture}
    \end{minipage}
  }
  \\
  \subcaptionbox{DDR}{%
    \begin{minipage}{0.48\textwidth}
\noindent\begin{tikzpicture}
\begin{loglogaxis}[
  width=\linewidth,
  height=0.78\linewidth,
  title={h1 potential norm vs. $h$},
  title style={font=\small},
  label style={font=\scriptsize},
  tick label style={font=\scriptsize},
  grid=both,
  major grid style={gray!40},
  minor grid style={gray!20},
  legend pos=south east,
  legend style={font=\scriptsize, inner xsep=2pt, inner ysep=1pt},
  legend cell align={left}
]
\addplot+[blue, mark=*, thick] coordinates { (6.6146214270788017e-1,4.8556813312870956e-1) (3.8885226009961682e-1,1.3548774667605548e-1) (2.0379771530539792e-1,5.4047195754479083e-2) (1.0586955040260637e-1,1.7536587987896465e-2) (5.5649832160483637e-2,5.3177213439192855e-3) };
\addlegendentry{$k = 0$}
\draw[draw=none] (axis cs:6.6146214270788017e-1,4.8556813312870956e-1) -- (axis cs:3.8885226009961682e-1,1.3548774667605548e-1) node[midway, above, text=blue, fill=white, fill opacity=0.8, text opacity=1, inner sep=1pt, font=\scriptsize] {1.84};
\draw[draw=none] (axis cs:3.8885226009961682e-1,1.3548774667605548e-1) -- (axis cs:2.0379771530539792e-1,5.4047195754479083e-2) node[midway, above, text=blue, fill=white, fill opacity=0.8, text opacity=1, inner sep=1pt, font=\scriptsize] {1.33};
\draw[draw=none] (axis cs:2.0379771530539792e-1,5.4047195754479083e-2) -- (axis cs:1.0586955040260637e-1,1.7536587987896465e-2) node[midway, above, text=blue, fill=white, fill opacity=0.8, text opacity=1, inner sep=1pt, font=\scriptsize] {1.62};
\draw[draw=none] (axis cs:1.0586955040260637e-1,1.7536587987896465e-2) -- (axis cs:5.5649832160483637e-2,5.3177213439192855e-3) node[midway, above, text=blue, fill=white, fill opacity=0.8, text opacity=1, inner sep=1pt, font=\scriptsize] {1.72};
\addplot+[red, mark=square*, thick] coordinates { (6.6146214270788017e-1,1.4379474915646581e-1) (3.8885226009961682e-1,3.4379803829323027e-2) (2.0379771530539792e-1,6.5860353577954145e-3) (1.0586955040260637e-1,1.1661361323514383e-3) (5.5649832160483637e-2,2.0753184140412910e-4) };
\addlegendentry{$k = 1$}
\draw[draw=none] (axis cs:6.6146214270788017e-1,1.4379474915646581e-1) -- (axis cs:3.8885226009961682e-1,3.4379803829323027e-2) node[midway, above, text=red, fill=white, fill opacity=0.8, text opacity=1, inner sep=1pt, font=\scriptsize] {2.06};
\draw[draw=none] (axis cs:3.8885226009961682e-1,3.4379803829323027e-2) -- (axis cs:2.0379771530539792e-1,6.5860353577954145e-3) node[midway, above, text=red, fill=white, fill opacity=0.8, text opacity=1, inner sep=1pt, font=\scriptsize] {2.38};
\draw[draw=none] (axis cs:2.0379771530539792e-1,6.5860353577954145e-3) -- (axis cs:1.0586955040260637e-1,1.1661361323514383e-3) node[midway, above, text=red, fill=white, fill opacity=0.8, text opacity=1, inner sep=1pt, font=\scriptsize] {2.50};
\draw[draw=none] (axis cs:1.0586955040260637e-1,1.1661361323514383e-3) -- (axis cs:5.5649832160483637e-2,2.0753184140412910e-4) node[midway, above, text=red, fill=white, fill opacity=0.8, text opacity=1, inner sep=1pt, font=\scriptsize] {2.49};
\addplot+[brown, mark=o, thick] coordinates { (6.6146214270788017e-1,3.1183427555162145e-2) (3.8885226009961682e-1,4.1769233468910280e-3) (2.0379771530539792e-1,6.0944078317174820e-4) (1.0586955040260637e-1,8.8311362154375636e-5) (5.5649832160483637e-2,1.2136772811316552e-5) };
\addlegendentry{$k = 2$}
\draw[draw=none] (axis cs:6.6146214270788017e-1,3.1183427555162145e-2) -- (axis cs:3.8885226009961682e-1,4.1769233468910280e-3) node[midway, above, text=brown, fill=white, fill opacity=0.8, text opacity=1, inner sep=1pt, font=\scriptsize] {2.90};
\draw[draw=none] (axis cs:3.8885226009961682e-1,4.1769233468910280e-3) -- (axis cs:2.0379771530539792e-1,6.0944078317174820e-4) node[midway, above, text=brown, fill=white, fill opacity=0.8, text opacity=1, inner sep=1pt, font=\scriptsize] {2.78};
\draw[draw=none] (axis cs:2.0379771530539792e-1,6.0944078317174820e-4) -- (axis cs:1.0586955040260637e-1,8.8311362154375636e-5) node[midway, above, text=brown, fill=white, fill opacity=0.8, text opacity=1, inner sep=1pt, font=\scriptsize] {2.79};
\draw[draw=none] (axis cs:1.0586955040260637e-1,8.8311362154375636e-5) -- (axis cs:5.5649832160483637e-2,1.2136772811316552e-5) node[midway, above, text=brown, fill=white, fill opacity=0.8, text opacity=1, inner sep=1pt, font=\scriptsize] {2.86};
\addplot+[teal, mark=triangle*, thick] coordinates { (6.6146214270788017e-1,3.7448360321259296e-3) (3.8885226009961682e-1,3.2466856078611076e-4) (2.0379771530539792e-1,2.0365806792985649e-5) (1.0586955040260637e-1,1.1514447597676716e-6) (5.5649832160483637e-2,6.5403514632437354e-8) };
\addlegendentry{$k = 3$}
\draw[draw=none] (axis cs:6.6146214270788017e-1,3.7448360321259296e-3) -- (axis cs:3.8885226009961682e-1,3.2466856078611076e-4) node[midway, above, text=teal, fill=white, fill opacity=0.8, text opacity=1, inner sep=1pt, font=\scriptsize] {3.53};
\draw[draw=none] (axis cs:3.8885226009961682e-1,3.2466856078611076e-4) -- (axis cs:2.0379771530539792e-1,2.0365806792985649e-5) node[midway, above, text=teal, fill=white, fill opacity=0.8, text opacity=1, inner sep=1pt, font=\scriptsize] {3.99};
\draw[draw=none] (axis cs:2.0379771530539792e-1,2.0365806792985649e-5) -- (axis cs:1.0586955040260637e-1,1.1514447597676716e-6) node[midway, above, text=teal, fill=white, fill opacity=0.8, text opacity=1, inner sep=1pt, font=\scriptsize] {4.14};
\draw[draw=none] (axis cs:1.0586955040260637e-1,1.1514447597676716e-6) -- (axis cs:5.5649832160483637e-2,6.5403514632437354e-8) node[midway, above, text=teal, fill=white, fill opacity=0.8, text opacity=1, inner sep=1pt, font=\scriptsize] {4.14};
\end{loglogaxis}
\end{tikzpicture}
    \end{minipage}
    \begin{minipage}{0.48\textwidth}
\noindent\begin{tikzpicture}
\begin{loglogaxis}[
  width=\linewidth,
  height=0.78\linewidth,
  title={l2 potential norm vs. $h$},
  title style={font=\small},
  label style={font=\scriptsize},
  tick label style={font=\scriptsize},
  grid=both,
  major grid style={gray!40},
  minor grid style={gray!20},
  legend pos=south east,
  legend style={font=\scriptsize, inner xsep=2pt, inner ysep=1pt},
  legend cell align={left}
]
\addplot+[blue, mark=*, thick] coordinates { (6.6146214270788017e-1,1.4815027292359575e-1) (3.8885226009961682e-1,2.4102211385402618e-2) (2.0379771530539792e-1,6.0888222966447823e-3) (1.0586955040260637e-1,1.4442636809378171e-3) (5.5649832160483637e-2,3.5385869969657355e-4) };
\addlegendentry{$k = 0$}
\draw[draw=none] (axis cs:6.6146214270788017e-1,1.4815027292359575e-1) -- (axis cs:3.8885226009961682e-1,2.4102211385402618e-2) node[midway, above, text=blue, fill=white, fill opacity=0.8, text opacity=1, inner sep=1pt, font=\scriptsize] {2.62};
\draw[draw=none] (axis cs:3.8885226009961682e-1,2.4102211385402618e-2) -- (axis cs:2.0379771530539792e-1,6.0888222966447823e-3) node[midway, above, text=blue, fill=white, fill opacity=0.8, text opacity=1, inner sep=1pt, font=\scriptsize] {1.98};
\draw[draw=none] (axis cs:2.0379771530539792e-1,6.0888222966447823e-3) -- (axis cs:1.0586955040260637e-1,1.4442636809378171e-3) node[midway, above, text=blue, fill=white, fill opacity=0.8, text opacity=1, inner sep=1pt, font=\scriptsize] {2.08};
\draw[draw=none] (axis cs:1.0586955040260637e-1,1.4442636809378171e-3) -- (axis cs:5.5649832160483637e-2,3.5385869969657355e-4) node[midway, above, text=blue, fill=white, fill opacity=0.8, text opacity=1, inner sep=1pt, font=\scriptsize] {2.03};
\addplot+[red, mark=square*, thick] coordinates { (6.6146214270788017e-1,5.4977360688685889e-2) (3.8885226009961682e-1,7.9572190048480193e-3) (2.0379771530539792e-1,8.2187522688792784e-4) (1.0586955040260637e-1,7.5724316024048654e-5) (5.5649832160483637e-2,7.2613375298368452e-6) };
\addlegendentry{$k = 1$}
\draw[draw=none] (axis cs:6.6146214270788017e-1,5.4977360688685889e-2) -- (axis cs:3.8885226009961682e-1,7.9572190048480193e-3) node[midway, above, text=red, fill=white, fill opacity=0.8, text opacity=1, inner sep=1pt, font=\scriptsize] {2.79};
\draw[draw=none] (axis cs:3.8885226009961682e-1,7.9572190048480193e-3) -- (axis cs:2.0379771530539792e-1,8.2187522688792784e-4) node[midway, above, text=red, fill=white, fill opacity=0.8, text opacity=1, inner sep=1pt, font=\scriptsize] {3.28};
\draw[draw=none] (axis cs:2.0379771530539792e-1,8.2187522688792784e-4) -- (axis cs:1.0586955040260637e-1,7.5724316024048654e-5) node[midway, above, text=red, fill=white, fill opacity=0.8, text opacity=1, inner sep=1pt, font=\scriptsize] {3.44};
\draw[draw=none] (axis cs:1.0586955040260637e-1,7.5724316024048654e-5) -- (axis cs:5.5649832160483637e-2,7.2613375298368452e-6) node[midway, above, text=red, fill=white, fill opacity=0.8, text opacity=1, inner sep=1pt, font=\scriptsize] {3.38};
\addplot+[brown, mark=o, thick] coordinates { (6.6146214270788017e-1,1.1861335558273463e-2) (3.8885226009961682e-1,8.6630961278745663e-4) (2.0379771530539792e-1,7.5339349178386771e-5) (1.0586955040260637e-1,6.2690505980241154e-6) (5.5649832160483637e-2,4.7202365984620344e-7) };
\addlegendentry{$k = 2$}
\draw[draw=none] (axis cs:6.6146214270788017e-1,1.1861335558273463e-2) -- (axis cs:3.8885226009961682e-1,8.6630961278745663e-4) node[midway, above, text=brown, fill=white, fill opacity=0.8, text opacity=1, inner sep=1pt, font=\scriptsize] {3.78};
\draw[draw=none] (axis cs:3.8885226009961682e-1,8.6630961278745663e-4) -- (axis cs:2.0379771530539792e-1,7.5339349178386771e-5) node[midway, above, text=brown, fill=white, fill opacity=0.8, text opacity=1, inner sep=1pt, font=\scriptsize] {3.52};
\draw[draw=none] (axis cs:2.0379771530539792e-1,7.5339349178386771e-5) -- (axis cs:1.0586955040260637e-1,6.2690505980241154e-6) node[midway, above, text=brown, fill=white, fill opacity=0.8, text opacity=1, inner sep=1pt, font=\scriptsize] {3.59};
\draw[draw=none] (axis cs:1.0586955040260637e-1,6.2690505980241154e-6) -- (axis cs:5.5649832160483637e-2,4.7202365984620344e-7) node[midway, above, text=brown, fill=white, fill opacity=0.8, text opacity=1, inner sep=1pt, font=\scriptsize] {3.73};
\addplot+[teal, mark=triangle*, thick] coordinates { (6.6146214270788017e-1,1.3587215825978463e-3) (3.8885226009961682e-1,6.4623909657233688e-5) (2.0379771530539792e-1,2.2060772426827392e-6) (1.0586955040260637e-1,6.0245333990790821e-8) (5.5649832160483637e-2,1.6523137553563396e-9) };
\addlegendentry{$k = 3$}
\draw[draw=none] (axis cs:6.6146214270788017e-1,1.3587215825978463e-3) -- (axis cs:3.8885226009961682e-1,6.4623909657233688e-5) node[midway, above, text=teal, fill=white, fill opacity=0.8, text opacity=1, inner sep=1pt, font=\scriptsize] {4.39};
\draw[draw=none] (axis cs:3.8885226009961682e-1,6.4623909657233688e-5) -- (axis cs:2.0379771530539792e-1,2.2060772426827392e-6) node[midway, above, text=teal, fill=white, fill opacity=0.8, text opacity=1, inner sep=1pt, font=\scriptsize] {4.87};
\draw[draw=none] (axis cs:2.0379771530539792e-1,2.2060772426827392e-6) -- (axis cs:1.0586955040260637e-1,6.0245333990790821e-8) node[midway, above, text=teal, fill=white, fill opacity=0.8, text opacity=1, inner sep=1pt, font=\scriptsize] {5.19};
\draw[draw=none] (axis cs:1.0586955040260637e-1,6.0245333990790821e-8) -- (axis cs:5.5649832160483637e-2,1.6523137553563396e-9) node[midway, above, text=teal, fill=white, fill opacity=0.8, text opacity=1, inner sep=1pt, font=\scriptsize] {5.19};
\end{loglogaxis}
\end{tikzpicture}
    \end{minipage}
  } 
  \caption{Comparison of HHO (above) and nodal DDR (below) solutions of the Poisson problem on a refined sequence of Voronoi meshes starting with the one depicted in Figure \ref{fig:voronoi}.}
  \label{fig:hho-ddr-poisson-comparison:voronoi}
\end{figure}

We first consider the solution of the Poisson problem using the HHO method and the variant of the nodal DDR method discussed in \cite[Section 3]{Beirao-da-Veiga.Di-Pietro.ea:26}.
Figures \ref{fig:hho-ddr-poisson-comparison:tria} and \ref{fig:hho-ddr-poisson-comparison:voronoi} show convergence on the same sequences of refined triangular and polygonal meshes for degrees from $0$ to $3$.
These figures were obtained using the \lstinline[style=GenericStyle]{--latex-plots} option of the \lstinline[style=GenericStyle]{convergence_test.sh} script discussed in Section \ref{sec:solvers.scripts:convergence_test.sh} below. The resulting plots were included without modifications in this document using the \lstinline[style=GenericStyle]{\input} \LaTeX command.
It is interesting to observe that, for such a simple problem, at every refinement level $i$, DDR is more accurate (and it can be checked that it also has fewer unknowns) than HHO, with superconvergence for the lowest degrees on the triangular mesh family.

\begin{figure}\centering
  \includegraphics[height=4.0cm]{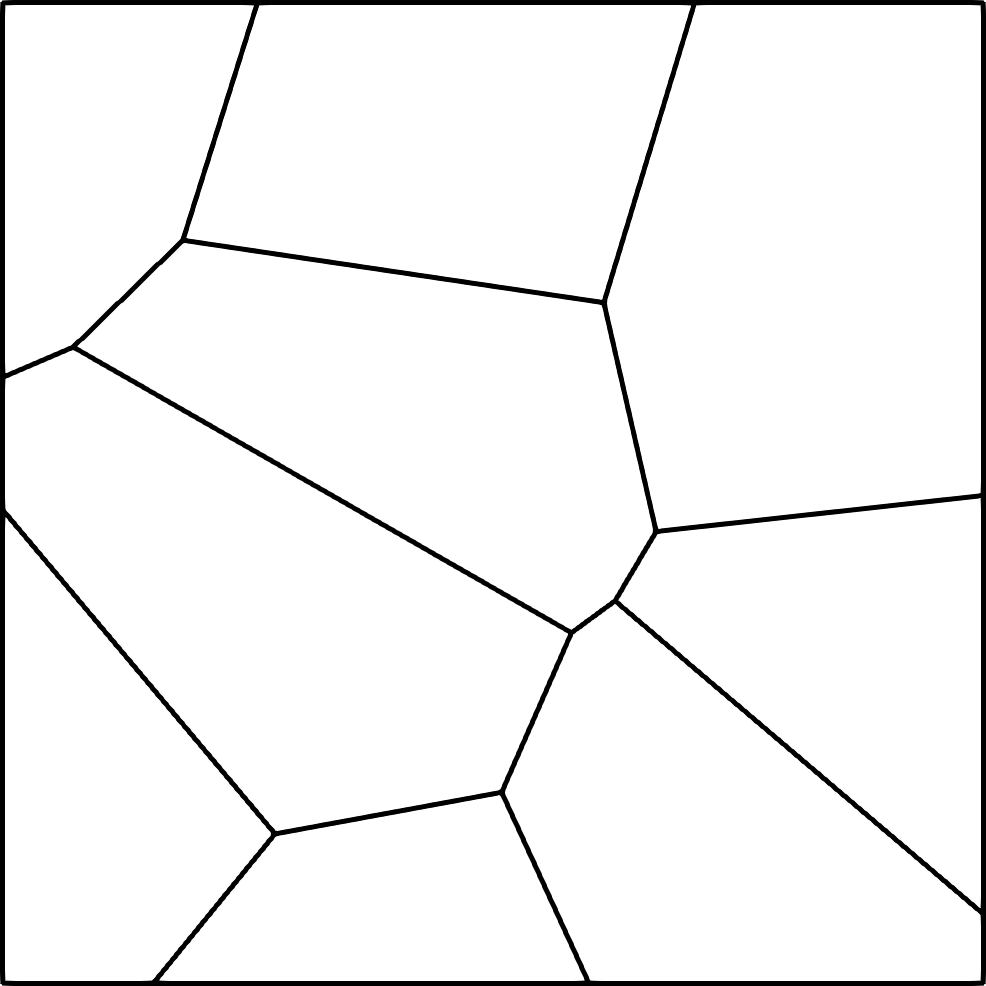}
  \caption{Polygonal mesh used for the degree-enrichment study of Figure \ref{fig:hho:degree.enrichment}.}
  \label{fig:voronoi}
\end{figure}

\begin{figure}\centering
  \subcaptionbox{$k=0$}{%
    \includegraphics[height=5.0cm]{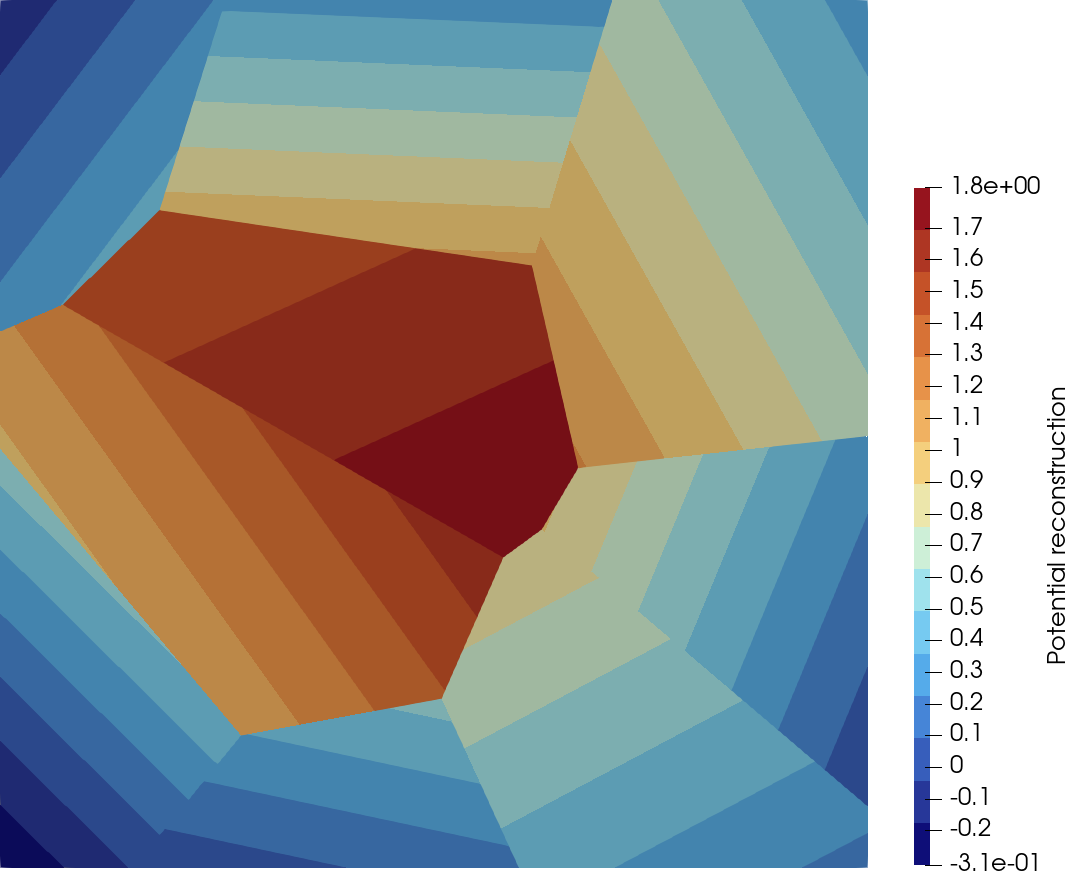}
  }
  \hspace{1cm}
  \subcaptionbox{$k=1$}{%
    \includegraphics[height=5.0cm]{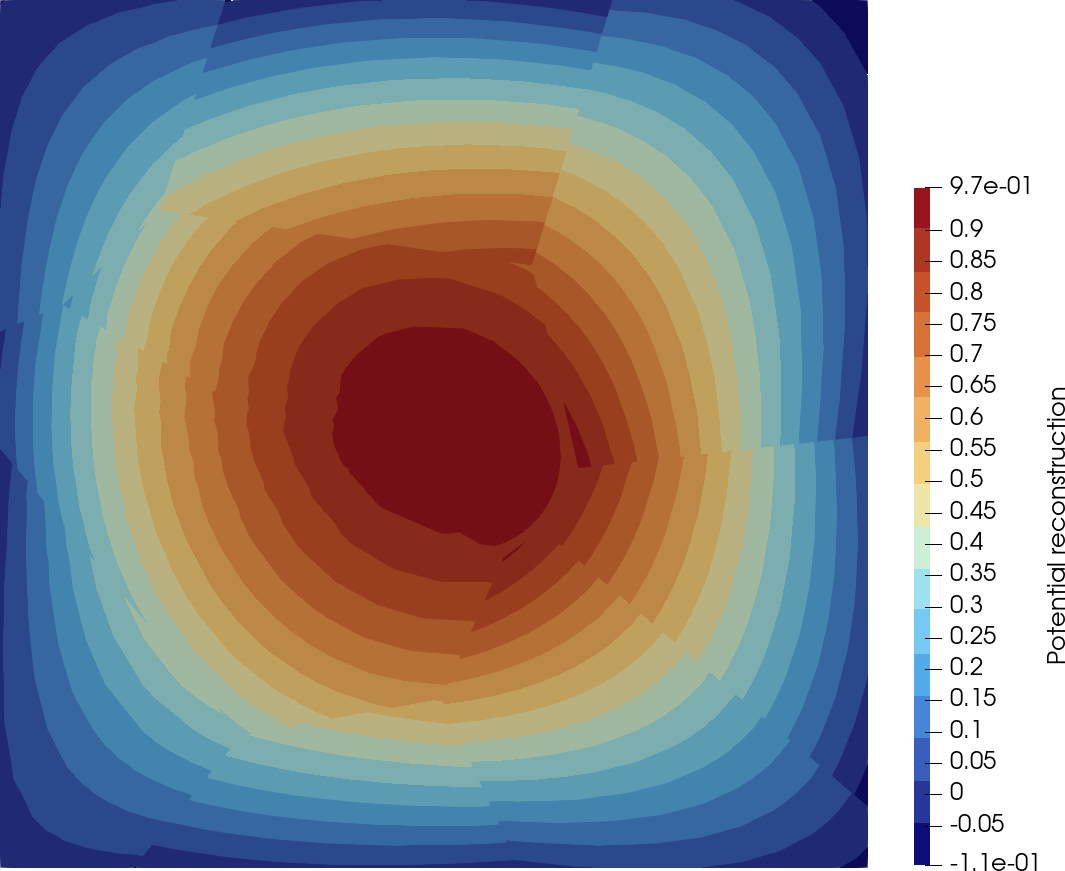}
  }
  \\[0.5cm]
  \subcaptionbox{$k=2$}{%
    \includegraphics[height=5.0cm]{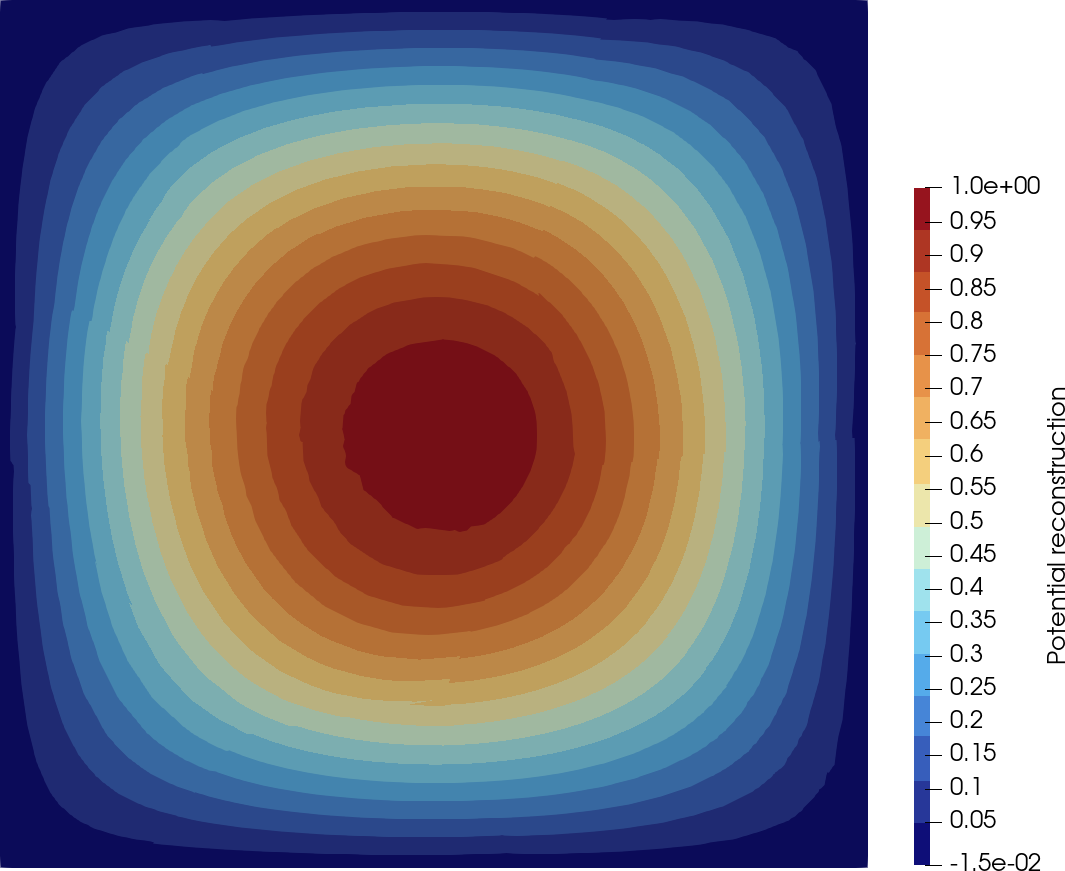}
  }
  \hspace{1cm}  
  \subcaptionbox{$k=3$}{%
    \includegraphics[height=5.0cm]{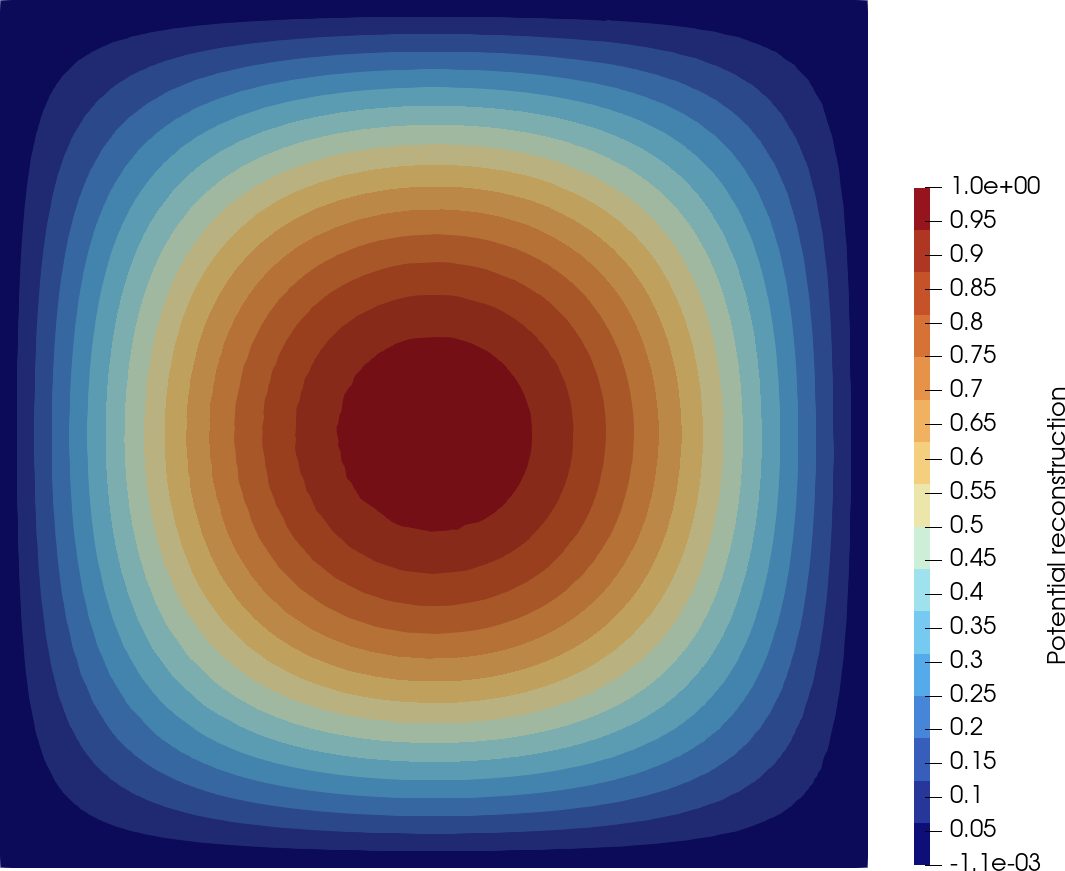}
  }
  \\[0.5cm]
  \subcaptionbox{$k=4$}{%
    \includegraphics[height=5.0cm]{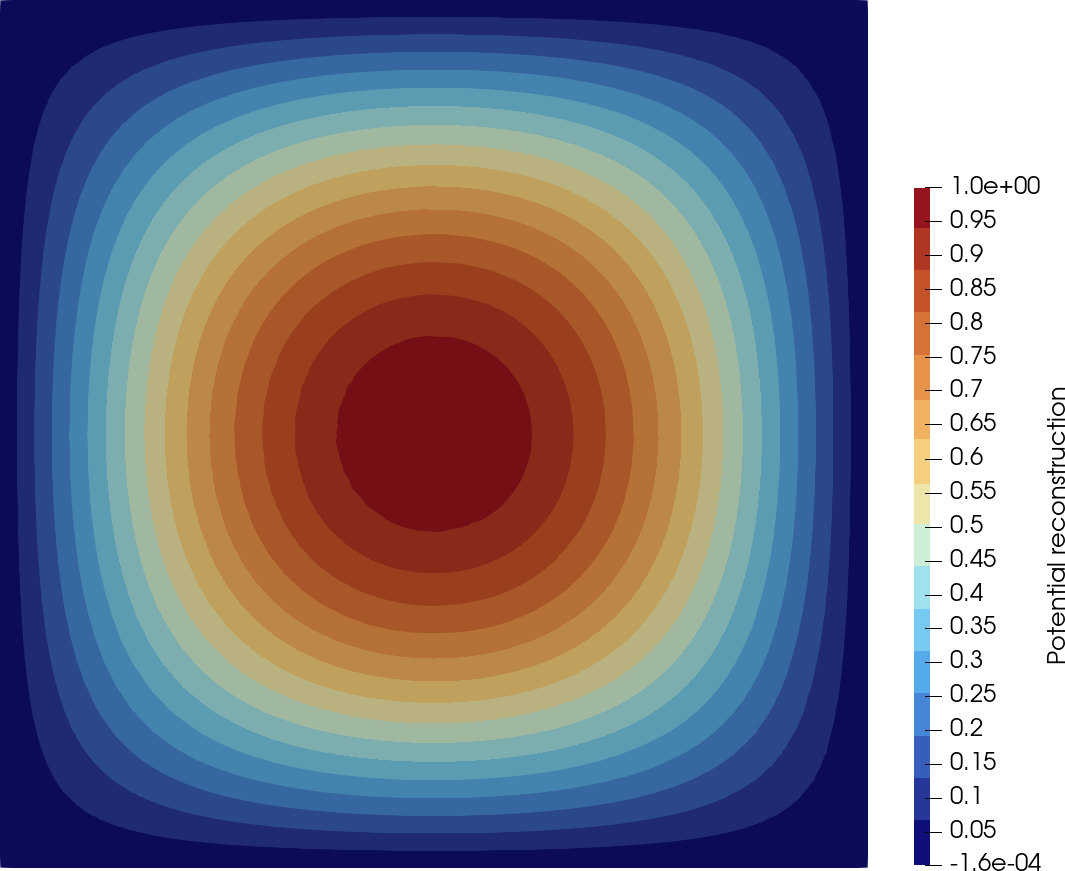}
  }
  \hspace{1cm}
  \subcaptionbox{$k=5$}{%
    \includegraphics[height=5.0cm]{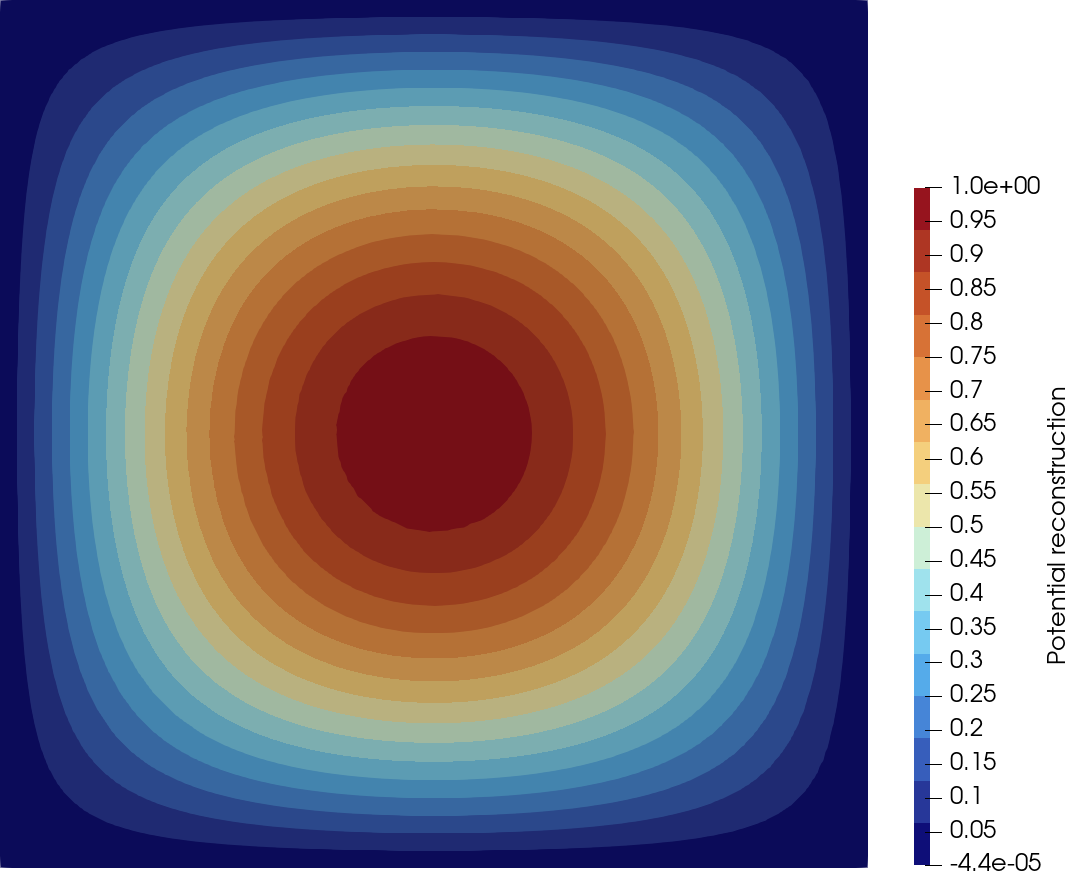}
  }
  \caption{Degree enrichment of the HHO method implemented in Section \ref{sec:overview} on the Voronoi mesh of Figure \ref{fig:voronoi}.}
  \label{fig:hho:degree.enrichment}
\end{figure}

The second example showcases the export features of \polyrust by displaying the HHO solution on the coarse polygonal mesh depicted in Figure \ref{fig:voronoi}.
Specifically, Figure \ref{fig:hho:degree.enrichment} contains a degree-enrichment study that displays the potential reconstruction for polynomial degrees ranging from $0$ up to $5$.
The export feature of \polyrust creates VTU files where the field is evaluated on a submesh to account for the variability of the high-degree solution inside each element.

\section{Solvers and scripts}\label{sec:solvers}

\polyrust offers a \Rust implementation of an interpreter for the DSL discussed in the previous section.
The interpreter's tools are intended for constructing solvers, i.e., \Rust executables that process and validate the DSL configuration file and use its facilities to simplify the implementation of complicated mathematical objects such as spaces and forms.
\polyrust comes with four basic problem-agnostic solvers: for steady/unsteady linear/nonlinear problems.
The focus here is on the steady solvers, which have been used to run all the numerical tests in this article.
We present them from the perspective of a basic user: the actual implementation of the DSL interpreter and the solvers is an advanced topic that lies outside the scope of the present article.
We also briefly describe a script for running convergence tests, which is particularly useful when prototyping numerical methods.

\subsection{Generic solver for linear problems}

The basic usage of the generic solver for linear problems from an installed version of \polyrust is
\begin{lstlisting}[style=GenericStyle, caption={Usage of the generic solver for linear problems}, label=lst:dsl_solver_linear_problem]
dsl_solver_linear_problem --dsl DSL_FILE --mesh MESH [OPTIONS]
\end{lstlisting}
Here, the mandatory arguments \lstinline[style=GenericStyle]{DSL_FILE} and \lstinline[style=GenericStyle]{MESH} are, respectively, the DSL configuration file and a VTK or VTU mesh file, while \lstinline[style=GenericStyle]{[OPTIONS]} is a list of command-line options described in Table \ref{tab:dsl-solver-linear-options}.
When running from compiled sources using \cargo, the invocation in Listing \ref{lst:dsl_solver_linear_problem} modifies as follows:
\begin{lstlisting}[style=GenericStyle]
cargo run --release --bin dsl_solver_linear_problem --features umfpack -- \
  --dsl DSL_FILE --mesh MESH [OPTIONS]
\end{lstlisting}
The \lstinline[style=GenericStyle]{--features umfpack} option activates the support for UMFPACK \cite{Davis:04}.
This feature is not enabled by default in order to maximize compatibility.
When the solver is compiled without it, the sparse algebraic problem is solved using a fallback sparse LU solver.

The solver processes a DSL file containing exactly one linear problem, which can be posed on either an atomic space or a Cartesian product space (see Section \ref{sec:advanced.features:cartesian.product.spaces} below). It constructs the operators, identifies the DOFs constrained by the strong boundary conditions, assembles the sparse matrix and right-hand side, and solves the resulting algebraic system. It then evaluates any error functionals requested through \texttt{compute errors} blocks and, when requested, exports the operators selected by the problem's \texttt{export} list. Operator exactness tests can instead be run in a diagnostic mode that does not assemble or solve the global problem; see Listing \ref{lst:operator.exactness}.
The main operations carried out by the solver (construction of the operators, assembly of the sparse system, evaluation of functionals) use thread-based parallelism.

\begin{table}
  \centering
  \caption{Command-line options for \texttt{dsl\_solver\_linear\_problem}. An asterisk marks a required option.}
  \label{tab:dsl-solver-linear-options}
  \small
  \begin{tabularx}{\textwidth}{
      @{}>{\raggedright\arraybackslash}p{0.38\textwidth}
      >{\raggedright\arraybackslash}X@{}}
    \toprule
    \textbf{Option} & \textbf{Description} \\
    \midrule

    \texttt{--dsl FILE}\textsuperscript{*}
    & Specifies the DSL configuration file. \\

    \texttt{--mesh FILE}\textsuperscript{*}
    & Specifies the VTK or VTU mesh file. \\

    \texttt{--degree \(k\)}, \texttt{-k \(k\)}
    & Sets the polynomial degree to \(k\). The default value is \(0\). \\

    \texttt{--output-prefix PREFIX},
    \texttt{-o PREFIX}
    & Sets the prefix of the exported VTU files. The default prefix is
    \texttt{dsl\_solution}. \\

    \texttt{--export-vtu}
    & Exports in VTU format the element-based operators named in the
    \texttt{export} list of the linear problem. \\

    \texttt{--quadrature-degree \(n\)}
    & Sets to \(n\) the quadrature degree used for non-polynomial
    contributions to the matrix and right-hand side. Polynomial
    contributions use fast polynomial integration. The default value is
    \(10\). \\

    \texttt{--functional-quadrature-degree \(n\)}
    & Sets to \(n\) the quadrature degree used for non-polynomial
    integrands in the error functionals. The default value is \(12\). \\

    \texttt{--test-operator-exactness}
    & Runs the operator exactness tests declared in the DSL without
    assembling or solving the global problem. \\

    \texttt{--help}, \texttt{-h}
    & Displays the command-line usage information and exits. \\

    \bottomrule
  \end{tabularx}
\end{table}

\subsection{Generic solver for nonlinear problems}

The basic usage of the generic solver for nonlinear problems is
\begin{lstlisting}[style=GenericStyle]
dsl_solver_nonlinear_problem --dsl DSL_FILE --mesh MESH [OPTIONS]
\end{lstlisting}
The mandatory \lstinline[style=GenericStyle]{--dsl} and \lstinline[style=GenericStyle]{--mesh} arguments have the same meanings as in the linear case. The remaining options are described in Table \ref{tab:dsl-solver-nonlinear-options}.
The fallback sparse LU solver and the thread-based parallelism mentioned for the linear solver are also used here.

The solver processes a DSL file containing exactly one nonlinear problem, which can be posed on either an atomic space or a Cartesian product space.
It constructs the operators, determines the DOFs constrained by the strong boundary conditions, and obtains an initial guess from the \texttt{initial\_bilinear} and \texttt{initial\_linear} forms when both are specified, or initializes the unconstrained DOFs to zero when they are absent. The interpolant of the exact solution can be selected instead. The nonlinear system is then solved by a damped hybrid strategy.
If a fixed-point form is available, up to five fixed-point iterations are first used as a warm-up; damped Newton iterations follow, with a line search that accepts only residual-decreasing steps.
If Newton fails or the iteration limit is reached, the solver returns to fixed-point iterations whenever the corresponding form is available. After convergence, it evaluates any error functionals requested through \texttt{compute errors} blocks and, when requested, exports the operators selected by the problem's \texttt{export} list.
The solver also provides diagnostics for checking the Jacobian and the residual at the interpolated exact solution, as well as the individual contributions to nonlinear forms for atomic problems.

\begin{table}
  \centering
  \caption{Command-line options for
    \texttt{dsl\_solver\_nonlinear\_problem}. An asterisk marks a required option.}
  \label{tab:dsl-solver-nonlinear-options}
  \small
  \begin{tabularx}{\textwidth}{
      @{}>{\raggedright\arraybackslash}p{0.43\textwidth}
      >{\raggedright\arraybackslash}X@{}}
    \toprule
    \textbf{Option} & \textbf{Description} \\
    \midrule

    \texttt{--dsl FILE}\textsuperscript{*}
    & Specifies the DSL configuration file. \\

    \texttt{--mesh FILE}\textsuperscript{*}
    & Specifies the VTK or VTU mesh file. \\

    \texttt{--degree \(k\)}, \texttt{-k \(k\)}
    & Sets the polynomial degree to \(k\). The default value is \(0\). \\

    \texttt{--initial-guess linear}
    & Uses the solution of the initial linear problem, when specified,
    as the initial guess. This is the default. \\

    \texttt{--initial-guess exact-interpolant}
    & Uses the interpolant associated with the error measures as the
    initial guess. \\

    \texttt{--output-prefix PREFIX},
    \texttt{-o PREFIX}
    & Sets the prefix of the exported VTU files. The default prefix is
    \texttt{dsl\_solution}. \\

    \texttt{--export-vtu}
    & Exports in VTU format the element-based operators named in the
    \texttt{export} list of the nonlinear problem. \\

    \texttt{--newton-tolerance TOL}
    & Sets the nonlinear residual tolerance. The default value is
    \(10^{-9}\). \\

    \texttt{--newton-max-iterations \(n\)}
    & Sets the maximum number of Newton iterations to \(n\). The default
    value is \(50\). \\

    \texttt{--nonlinear-quadrature-degree \(n\)}
    & Sets to \(n\) the quadrature degree used for non-polynomial
    contributions to nonlinear residuals and Jacobians. Polynomial
    contributions use fast polynomial integration. The default value is
    \(10\). \\

    \texttt{--functional-quadrature-degree \(n\)}
    & Sets to \(n\) the quadrature degree used for non-polynomial
    integrands in the error functionals. The default value is \(12\). \\

    \texttt{--check-jacobian}
    & Compares the Jacobian at the interpolated exact solution with
    centered finite differences in a deterministic direction. \\

    \texttt{--check-interpolant-residual}
    & Reports residual diagnostics at the interpolated exact solution. \\

    \texttt{--diagnose-nonlinear-contributions}
    & For an atomic problem, reports term-by-term residual and Jacobian
    diagnostics at the interpolated exact solution. \\

    \texttt{--test-operator-exactness}
    & Runs the operator exactness tests declared in the DSL without
    assembling or solving the nonlinear problem. \\

    \texttt{--help}, \texttt{-h}
    & Displays the command-line usage information and exits. \\

    \bottomrule
  \end{tabularx}
\end{table}

\subsection{Script for convergence tests}\label{sec:solvers.scripts:convergence_test.sh}

\polyrust also offers the \texttt{convergence\_test.sh} convenience script to automate convergence studies on refined mesh sequences.
Its basic usage is
\begin{lstlisting}[style=GenericStyle]
convergence_test.sh --dsl DSL_FILE --mesh-family MESH_FAMILY [OPTIONS]
\end{lstlisting}
The mandatory arguments \lstinline[style=GenericStyle]{DSL_FILE} and \lstinline[style=GenericStyle]{MESH_FAMILY} respectively correspond to a DSL description file and to the root name of a family of refined meshes.
For each polynomial degree and refinement level in user-selectable ranges, the script locates the corresponding VTK or VTU mesh family assuming the names \lstinline[style=GenericStyle]{MESH_FAMILY_i.vtu} or \lstinline[style=GenericStyle]{MESH_FAMILY_i.vtk}, and extracts the mesh size, the algebraic problem size, and all the error measures reported by the selected solver. In a source checkout, the solver is run through \cargo in release mode with UMFPACK enabled; in an installed copy, the selected executable is invoked directly.
It then computes the experimental convergence rates between consecutive refinement levels and collects the results in aligned tables.
The command-line options are listed in Table \ref{tab:convergence-test-options}.

The script always writes a text report. It can additionally produce a single \LaTeX{} table, either as a fragment suitable for inclusion using \lstinline[style=GenericStyle]{\input} or as a standalone document.
This is precisely how Tables \ref{tab:kovasznay:hho} and \ref{tab:kovasznay:hypre} were  generated.
It can also use \texttt{pgfplots} to generate log--log convergence plots such as the ones in Figures \ref{fig:hho-ddr-poisson-comparison:tria} and \ref{fig:hho-ddr-poisson-comparison:voronoi}.
The \lstinline[style=GenericStyle]{--latex-plots} option creates one \LaTeX{} fragment per error measure.
The \lstinline[style=GenericStyle]{--latex-plots-document} option instead groups all plots in a standalone, compilable document.
All output files are written to the current directory with names derived from the basename of the DSL configuration file and the selected solver.

\begin{table}
  \centering
  \caption{Command-line options for \texttt{convergence\_test.sh}. An asterisk marks a required option.}
  \label{tab:convergence-test-options}
  \small
  \begin{tabularx}{\textwidth}{
      @{}>{\raggedright\arraybackslash}p{0.40\textwidth}
      >{\raggedright\arraybackslash}X@{}}
    \toprule
    \textbf{Option} & \textbf{Description} \\
    \midrule

    \texttt{--dsl FILE}\textsuperscript{*}
    & Specifies the DSL configuration file. \\

    \texttt{--mesh-family NAME}\textsuperscript{*}
    & Selects meshes named \texttt{meshes/NAME\_i.vtk} or
    \texttt{meshes/NAME\_i.vtu}; the extension is detected
    automatically. \\

    \texttt{--bin FILE}
    & Selects the solver executable. The default is
    \texttt{dsl\_solver\_linear\_problem}. \\

    \texttt{--min-degree \(n\)},
    \texttt{--max-degree \(n\)}
    & Set the inclusive range of polynomial degrees. The default range
    is \(0,\ldots,2\). \\

    \texttt{--min-refinement-level \(i\)},
    \texttt{--max-refinement-level \(i\)}
    & Set the inclusive range of mesh refinement levels. The default
    range is \(0,\ldots,3\). \\

    \texttt{--args STRING}
    & Passes additional whitespace-separated arguments to the selected
    solver. \\

    \texttt{--export-vtu}
    & Exports the solution fields only for the maximum degree on the
    finest mesh. \\

    \texttt{--latex-tables}
    & Generates a \LaTeX{} fragment containing one convergence table for
    all polynomial degrees. \\

    \texttt{--latex-tables-document}
    & Generates a standalone, compilable \LaTeX{} document containing
    the convergence table. \\

    \texttt{--latex-plots}
    & Generates one \texttt{pgfplots} fragment per error measure, with
    all degrees and segment-wise estimated rates. \\

    \texttt{--latex-plots-document},
    \texttt{--latex-plot-document}
    & Generates a standalone, compilable \LaTeX{} document containing
    all convergence plots. The singular spelling is an alias. \\

    \texttt{--help}, \texttt{-h}
    & Displays the command-line usage information and exits. \\

    \bottomrule
  \end{tabularx}
\end{table}


\section{Advanced features}\label{sec:advanced.features}

In this section we discuss some advanced features of the language through concrete examples.
The complete DSL configuration files used for the examples in this section are available from the \polyrust website at \url{https://polyrust.org}.

\subsection{Cartesian product spaces}\label{sec:advanced.features:cartesian.product.spaces}

Cartesian product spaces are useful for multi-field problems and can be defined in the DSL using the syntax \lstinline{product space SPACE1 times ... times SPACEn}.
The following listing shows how to construct the Cartesian product space for the HHO discretization of the Stokes problem (see, e.g., \cite{Di-Pietro.Ern.ea:16*1} and also \cite[Chapter 8]{Di-Pietro.Droniou:20}):
\begin{lstlisting}
space Uh {
  element Poly(k, vector)
  edge Poly(k, vector)
}

space Ph {
  element Poly(k, scalar)
}

product space Xh = Uh times Ph
\end{lstlisting}
Bilinear forms can then be defined on any product of two Cartesian factors.
For example, the velocity-pressure coupling bilinear form for the above-mentioned method reads
\begin{lstlisting}[caption={Velocity-pressure coupling bilinear form}, label=lst:velocity.pressure]
bilinear form velocity_pressure : Uh(trial u) times Ph(test q) {
  sum_elements(-int(T) dof(u, T) dot grad(q) + int(dT) (dof(u, E) dot normal) * q)
}
\end{lstlisting}
and corresponds to rectangular local matrices.

The possibility of having multiple named DOFs supported by the same type of mesh entity also makes it possible to implement the method using a single space:
\begin{lstlisting}
space Xh {
  domain Poly(0, scalar)
  element Poly(k, vector) called element_velocity
  element Poly(k, scalar) called element_pressure
  edge Poly(k, vector)
}\end{lstlisting}
In this case, the velocity-pressure coupling bilinear form reads
\begin{lstlisting}
bilinear form velocity_pressure : Xh(trial u) times Xh(test v) {
  sum_elements(
               -int(T) dof(u, T, element_velocity) dot grad(dof(v, T, element_pressure))
               + int(dT) (dof(u, E) dot normal) * dof(v, T, element_pressure)
               )
}
\end{lstlisting}
The main difference with respect to the bilinear form in Listing \ref{lst:velocity.pressure} is that \lstinline{velocity_pressure} now results in larger square local matrices with zero blocks and hence in a potential efficiency penalty.
The implementation of the assembly step in \polyrust strives to minimize this penalty by first assembling the non-zero rectangular block and then extending it to a square matrix.
This process hinges on the identification of active DOFs based on the expression of the bilinear form and is automatic.

\begin{table}
  \centering
  \begin{tabular}{@{}r r r r r r r r r@{}}
    \toprule
    $i$ & $h$ & size & \multicolumn{2}{c}{velocity h1 norm} & \multicolumn{2}{c}{velocity l2 norm} & \multicolumn{2}{c}{pressure l2 norm} \\
    \cmidrule(lr){4-5}
    \cmidrule(lr){6-7}
    \cmidrule(lr){8-9}
    & &  & error & rate & error & rate & error & rate \\
    \midrule
    \multicolumn{9}{c}{$k = 0$} \\
    1 & 3.5355e-01 & 736 & 4.3351e+00 & -- & 1.6481e+00 & -- & 1.5461e-01 & -- \\
    2 & 1.7678e-01 & 3008 & 2.9067e+00 & 0.58 & 5.5621e-01 & 1.57 & 7.1725e-02 & 1.11 \\
    3 & 8.8388e-02 & 12160 & 1.5853e+00 & 0.87 & 1.5297e-01 & 1.86 & 3.0772e-02 & 1.22 \\
    4 & 4.4194e-02 & 48896 & 7.9540e-01 & 1.00 & 3.8816e-02 & 1.98 & 1.3839e-02 & 1.15 \\
    \midrule
    \multicolumn{9}{c}{$k = 1$} \\
    1 & 3.5355e-01 & 1856 & 1.7135e+00 & -- & 3.1462e-01 & -- & 2.7420e-02 & -- \\
    2 & 1.7678e-01 & 7552 & 4.7569e-01 & 1.85 & 4.6311e-02 & 2.76 & 5.9505e-03 & 2.20 \\
    3 & 8.8388e-02 & 30464 & 1.1923e-01 & 2.00 & 5.8691e-03 & 2.98 & 1.4300e-03 & 2.06 \\
    4 & 4.4194e-02 & 122368 & 2.8157e-02 & 2.08 & 6.8695e-04 & 3.09 & 3.7818e-04 & 1.92 \\
    \bottomrule
  \end{tabular}%
  \caption{Convergence results for the HHO method using the velocity-pressure coupling bilinear form defined in Listing \ref{lst:velocity.pressure} applied to the Kovasznay problem.}
  \label{tab:kovasznay:hho}
\end{table}

\begin{figure}\centering
  \includegraphics[height=6cm]{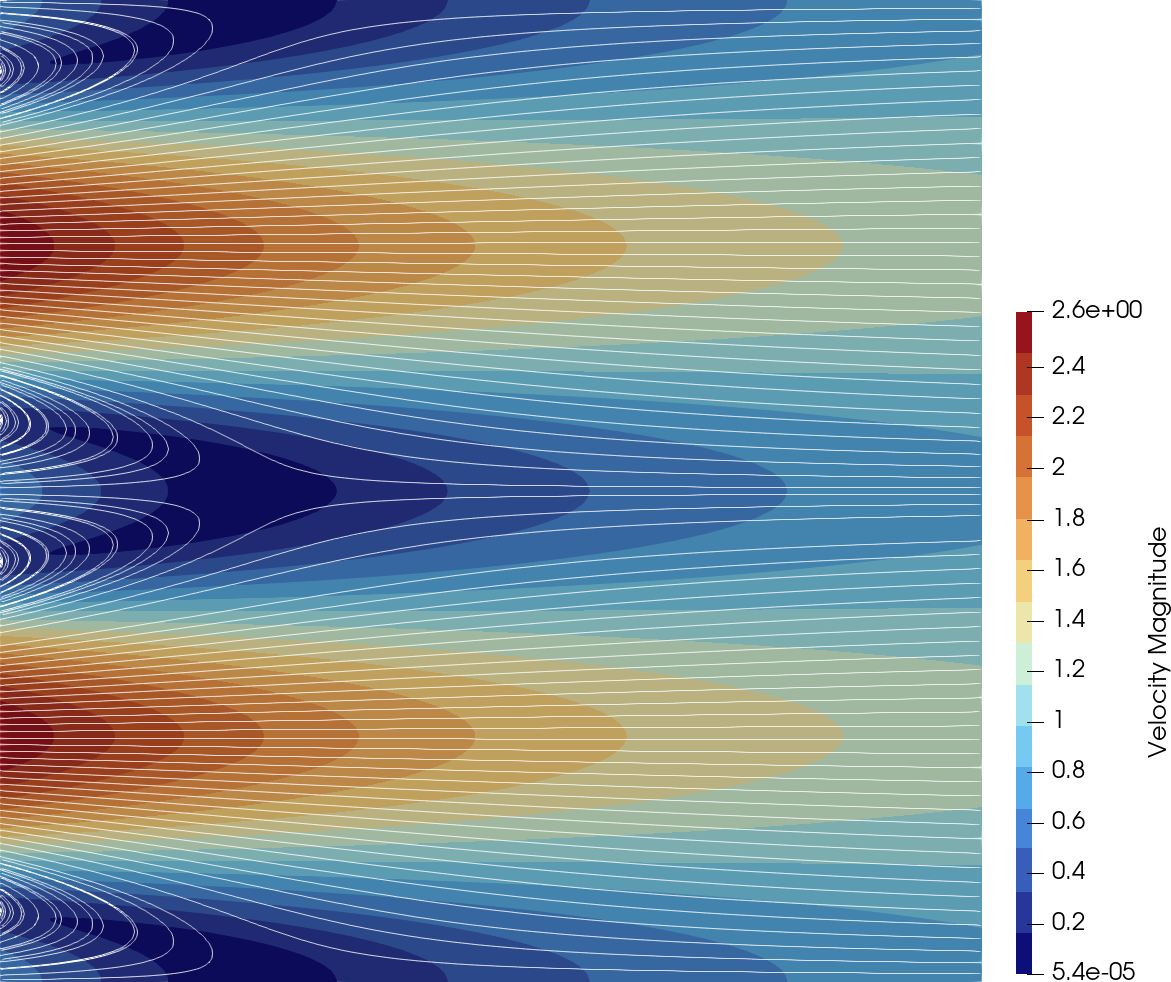}
  \hspace{1cm}
  \includegraphics[height=6cm]{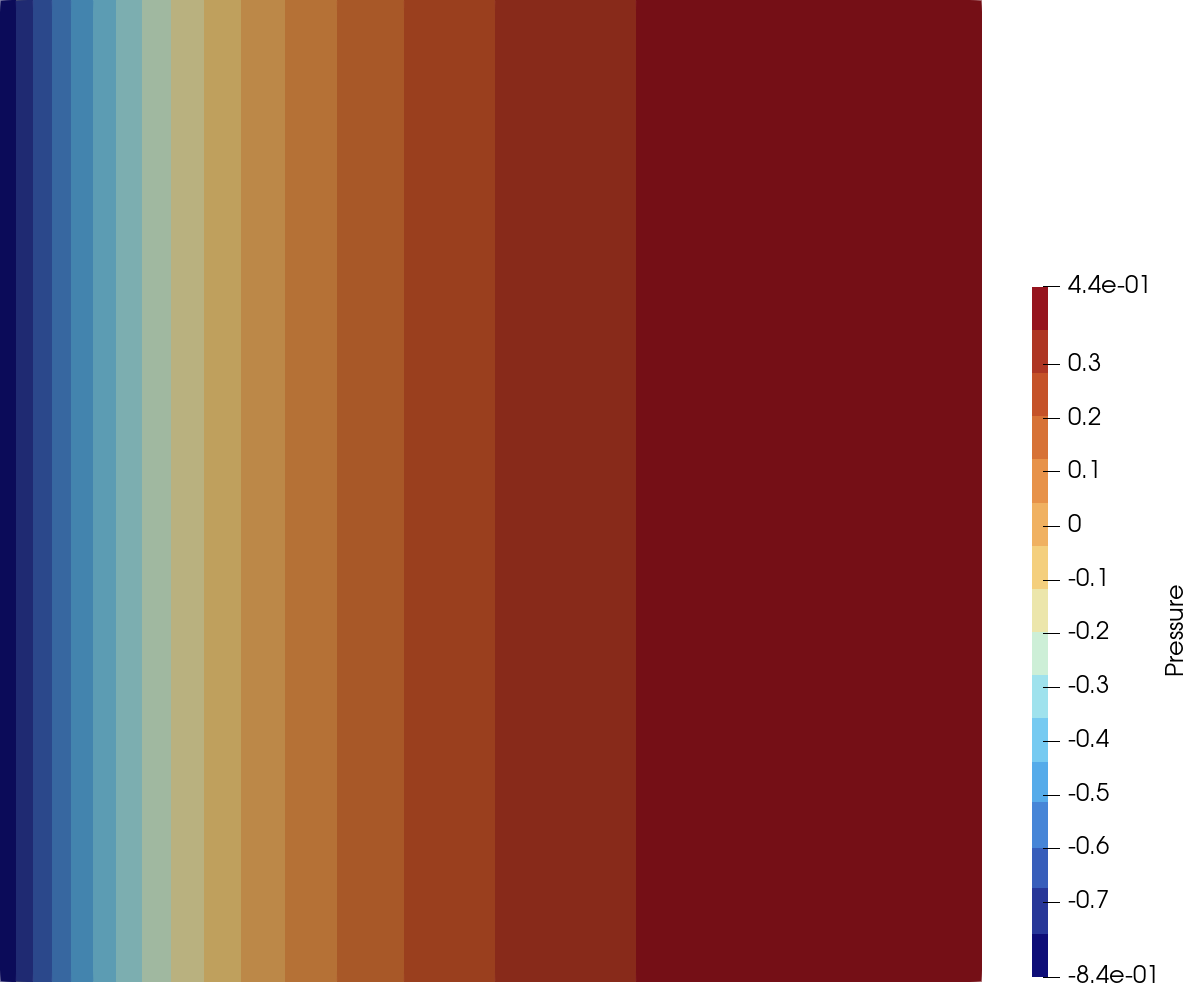}
  \caption{Velocity streamlines and magnitude (left) and pressure magnitude (right) for the Kovasznay solution used in the tests of Tables \ref{tab:kovasznay:hho} and \ref{tab:kovasznay:hypre}.}
  \label{fig:kovasznay}
\end{figure}

In Table \ref{tab:kovasznay:hho} we report convergence results for the HHO method of \cite{Botti.Di-Pietro.ea:19*1} with strongly enforced boundary conditions applied to the Kovasznay problem \cite{Kovasznay:48}; see Figure \ref{fig:kovasznay}.
The table was obtained using the convenience convergence script detailed in Section \ref{sec:solvers.scripts:convergence_test.sh}.

\subsection{Edge operators with vertex terms}

Virtual Element and DDR methods often require the reconstruction of operators on edges by mimicking integration by parts formulas.
These operators involve tangential derivatives, vertex evaluations, and possibly multiplication by the tangent or normal vectors to the edge.
The following listing shows an example of a complex edge operator resulting from the construction of \cite{Di-Pietro.Droniou.ea:25}:
\begin{lstlisting}
operator edge_hessian : Uh(u) -> Poly(k+1, vector) on edge E {
  forall w in Poly(k+1, vector):
    int(E) edge_hessian(u) dot w =
      int(E) dof(u, E, edge_value) * (tangential_derivative(tangential_derivative(w)) dot tangent(E))
      - int(E) dof(u, E, edge_normal_derivative) * (tangential_derivative(w) dot normal(E))
      + sum_vertices(orientation(V, E) * (
                                          dof(u, V, vertex_gradient) dot w(V)
                                          - dof(u, V, vertex_value) * (tangential_derivative(w) dot tangent(E))
                                         )
                    )
}
\end{lstlisting}

\begin{figure}\centering
  \subcaptionbox{HHO, $i = 2$}{%
    \includegraphics[height=5.0cm]{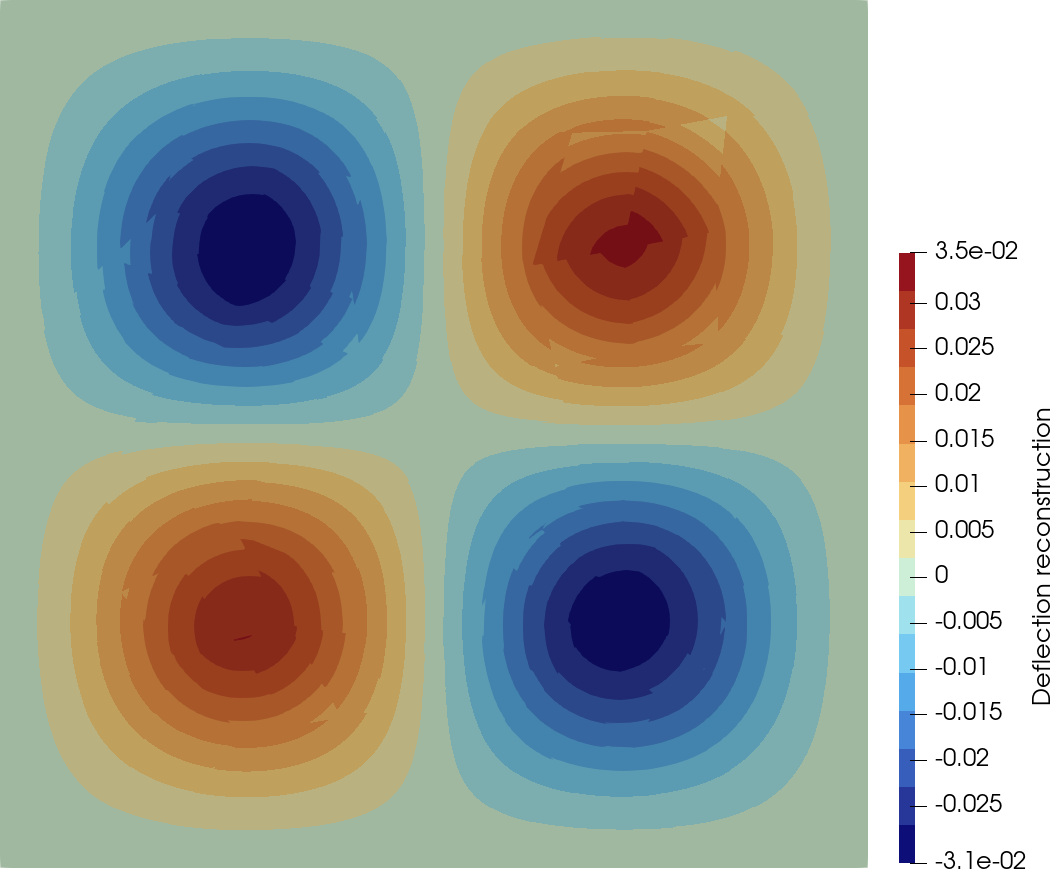}
  }
  \hspace{1cm}
  \subcaptionbox{DDR, $i = 2$}{%
    \includegraphics[height=5.0cm]{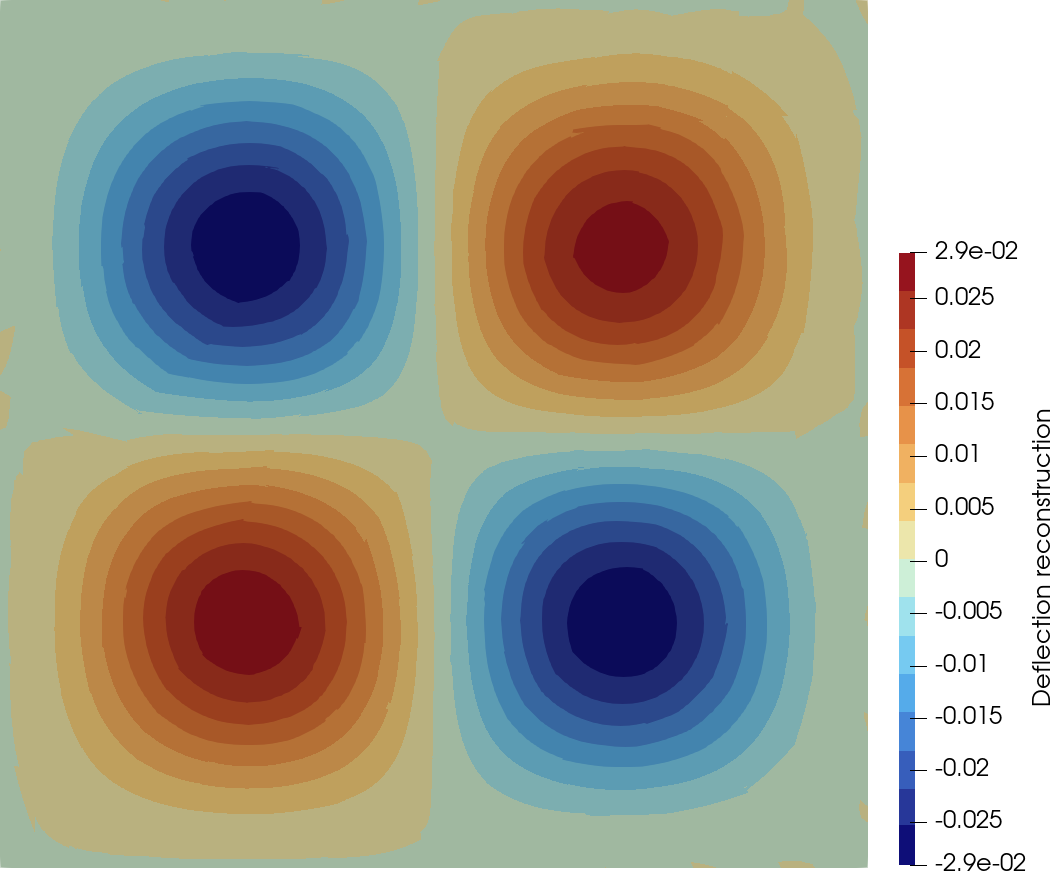}
  }
  \\[0.5cm]
  \subcaptionbox{HHO, $i = 3$}{%
    \includegraphics[height=5.0cm]{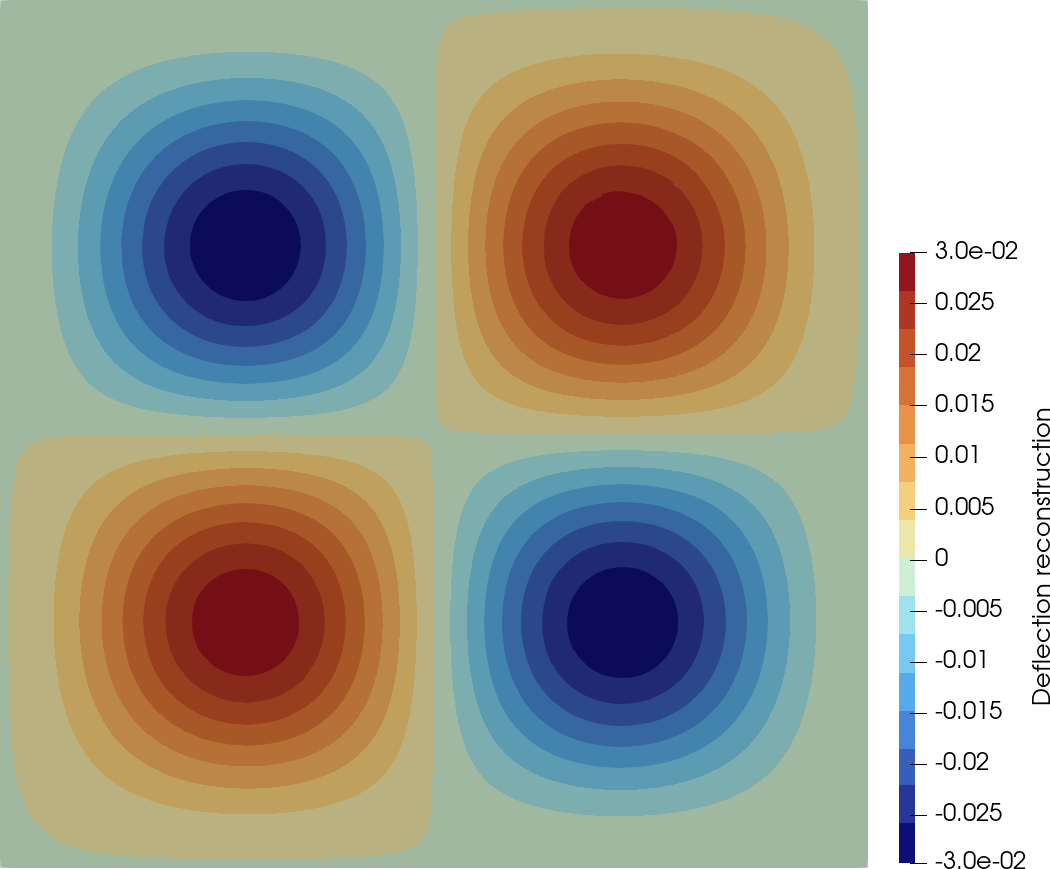}
  }
  \hspace{1cm}
  \subcaptionbox{DDR, $i = 3$}{%
    \includegraphics[height=5.0cm]{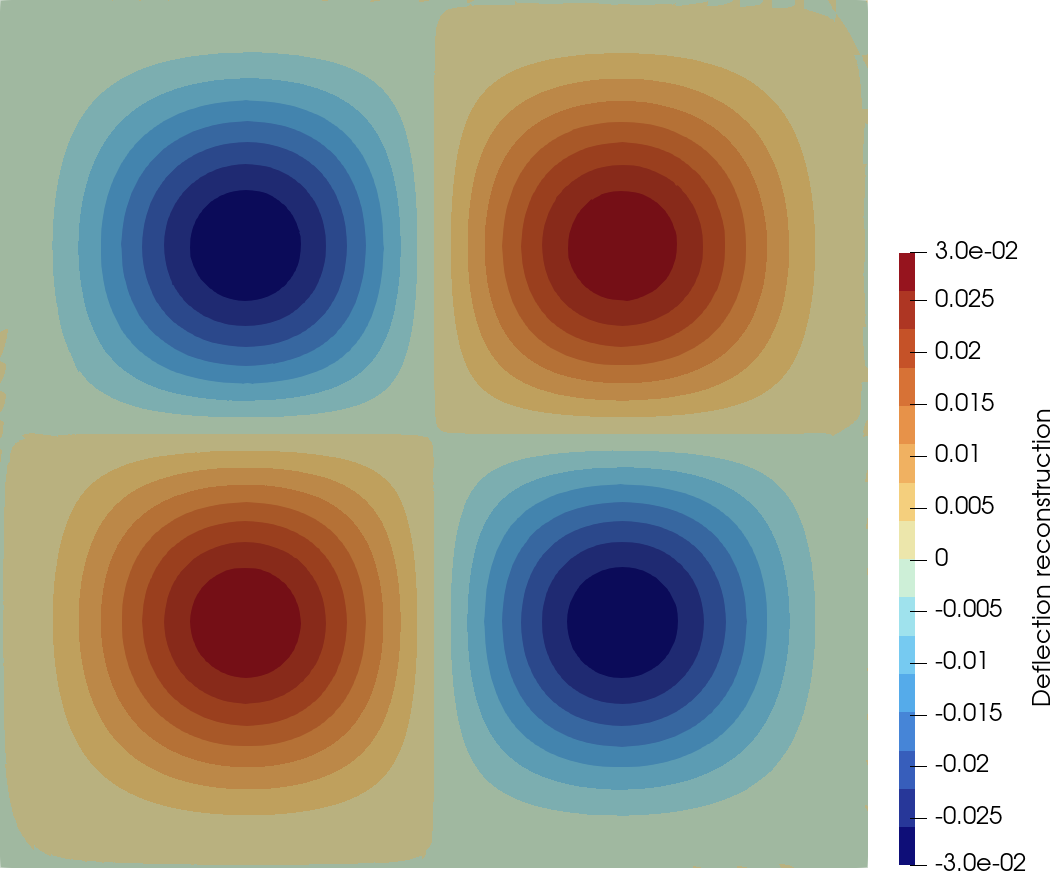}
  }
  \\[0.5cm]
  \subcaptionbox{HHO, $i = 4$}{%
    \includegraphics[height=5.0cm]{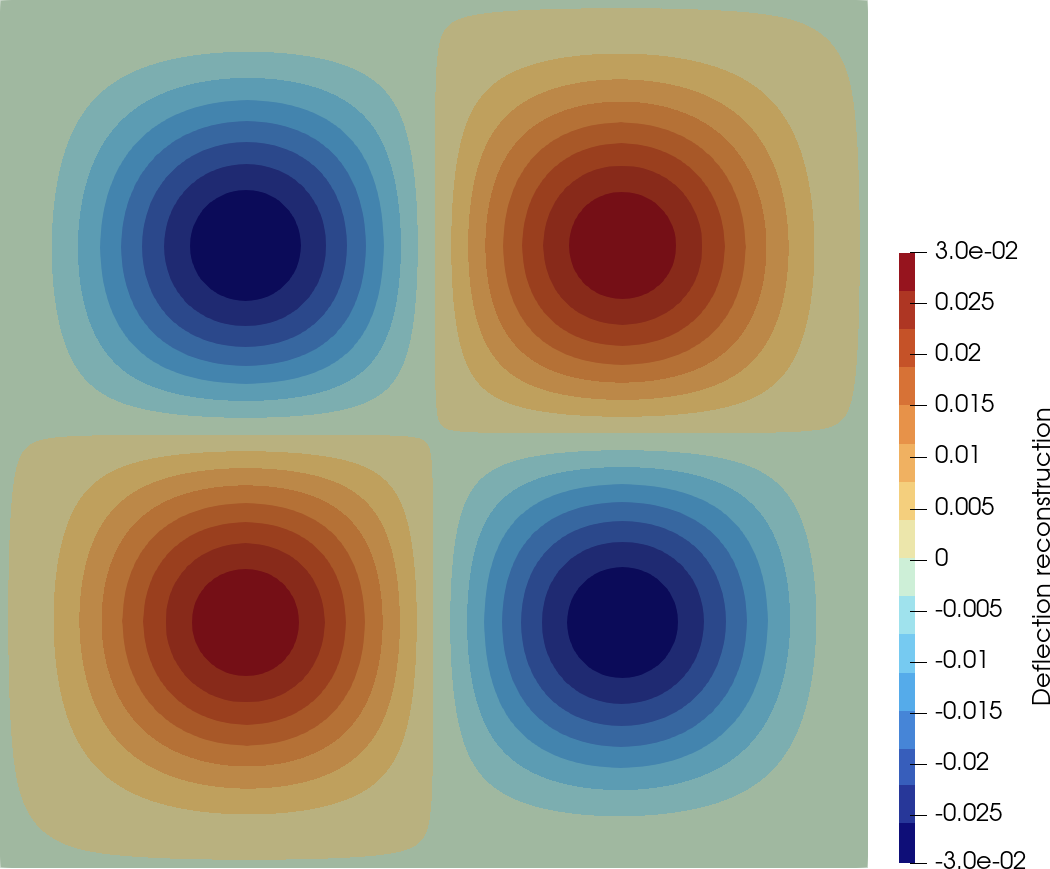}
  }
  \hspace{1cm}
  \subcaptionbox{DDR, $i = 4$}{%
    \includegraphics[height=5.0cm]{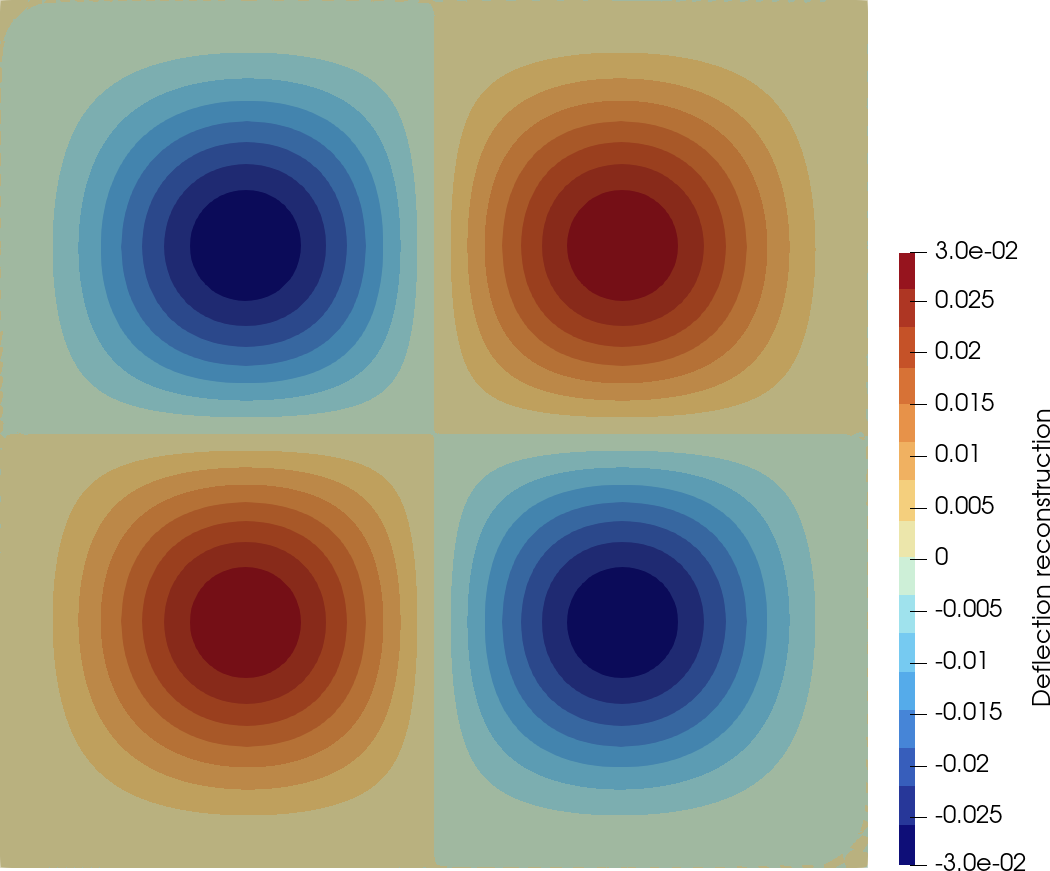}
  }
  \caption{Comparison of the HHO (left) and DDR (right) methods for the four-point load Kirchhoff--Love problem on a sequence of meshes obtained by refining the one in Figure \ref{fig:voronoi}.}
  \label{fig:hho-ddr-kirchhoff-love}
\end{figure}

In Figure \ref{fig:hho-ddr-kirchhoff-love} we compare the HHO method for Kirchhoff--Love plates of \cite{Bonaldi.Di-Pietro.ea:18} with polynomial degree $k = 2$ and a variation of the method of \cite{Di-Pietro.Droniou.ea:25} with polynomial degree $k = 0$ and higher-order discretization of the right-hand side.
With this choice of polynomial degrees, the methods have the same expected order of convergence in the energy norm.
The test corresponds to a four-point load configuration and is run on successive refinements of the mesh depicted in Figure \ref{fig:voronoi} (which corresponds to the refinement level $i = 0$).
The richer space structure of the DDR method appears to deliver better results on the coarsest mesh, whereas the two numerical solutions are visually indistinguishable on the finer mesh.

\subsection{Coefficient-dependent operators}

Operators can depend on spatially varying coefficients.
An example is the convective derivative reconstruction defined in \cite[Eq. (12)]{Di-Pietro.Droniou.ea:15}, which can be obtained in the DSL through the following definition:
\begin{lstlisting}
function velocity(vector X) -> vector = vector(0.5 - X[1], X[0] - 0.5)
  
operator advective_derivative : Uh(u) -> Poly(k, scalar) on element T {
  forall q in Poly(k, scalar):
    int(T) advective_derivative(u) * q =
    -int(T) dof(u, T) * (velocity dot grad(q))
    + int(dT) (velocity dot normal) * dof(u, E) * q
}
\end{lstlisting}

The integration kernel in \polyrust implements both standard quadratures from the \texttt{quadraturerules} crate and fast polynomial integration \cite{Chin.Lasserre.ea:15}.
The kernel automatically detects whether the integrand is a polynomial (including, e.g., the case where it is the product of polynomials and constant coefficients and when products of linear combinations of polynomials are present) and uses fast integration in this case.
Quadrature rules are used only for non-polynomial integrands, as in the \lstinline{advective_derivative} operator above or the \lstinline{load} linear form implemented in Listing \ref{lst:load}.
    
\subsection{Operators between polynomial spaces}

A syntax similar to that introduced in Section \ref{sec:overview:operators} can be used to define operators between polynomial spaces.
An example in which operators of this kind are needed is provided by the HYPRE method of \cite[Section 5.4]{Botti.Botti.ea:26} (Stokes) and \cite{Beirao-da-Veiga.Di-Pietro.ea:25} (Navier--Stokes), where the penalty term requires applying the Raviart--Thomas interpolator to vector-valued polynomials of degree up to $k+1$.
Its definition in the DSL is detailed in the following listing:
\begin{lstlisting}[caption={Raviart--Thomas interpolator in $\Poly{k+1}(T; \Real^2)$}, label=lst:raviart.thomas.interpolator]
operator raviart_thomas_interpolator_polykpo
  : Poly(k+1, vector)(u) -> RaviartThomasPoly(k+1) on element T {
  forall v in Poly(k-1, vector):
    int(T) raviart_thomas_interpolator_polykpo(u) dot v = int(T) u dot v
  forall edge E, forall q in Poly(k, scalar):
    int(E) (raviart_thomas_interpolator_polykpo(u) dot normal) * q
    = int(E) (u dot normal) * q
}
\end{lstlisting}
This listing incidentally shows another feature of operators on elements, namely the possibility to express variational conditions on the element edges.

\begin{table}
  \centering
  \begin{tabular}{@{}r r r r r r r r r@{}}
    \toprule
    $i$ & $h$ & size & \multicolumn{2}{c}{velocity h1 norm} & \multicolumn{2}{c}{velocity l2 norm} & \multicolumn{2}{c}{pressure l2 norm} \\
    \cmidrule(lr){4-5}
    \cmidrule(lr){6-7}
    \cmidrule(lr){8-9}
    & &  & error & rate & error & rate & error & rate \\
    \midrule
    \multicolumn{9}{c}{$k = 0$} \\
    1 & 3.5355e-01 & 1073 & 2.7186e+00 & -- & 9.7862e-01 & -- & 1.8236e-01 & -- \\
    2 & 1.7678e-01 & 4321 & 2.1961e+00 & 0.31 & 4.0189e-01 & 1.28 & 9.1591e-02 & 0.99 \\
    3 & 8.8388e-02 & 17345 & 1.4356e+00 & 0.61 & 1.3301e-01 & 1.60 & 4.1867e-02 & 1.13 \\
    4 & 4.4194e-02 & 69505 & 8.5016e-01 & 0.76 & 4.0344e-02 & 1.72 & 1.9450e-02 & 1.11 \\
    \midrule
    \multicolumn{9}{c}{$k = 1$} \\
    1 & 3.5355e-01 & 2529 & 9.8806e-01 & -- & 1.6405e-01 & -- & 3.0343e-02 & -- \\
    2 & 1.7678e-01 & 10177 & 2.3755e-01 & 2.06 & 2.2946e-02 & 2.84 & 5.7044e-03 & 2.41 \\
    3 & 8.8388e-02 & 40833 & 6.6929e-02 & 1.83 & 3.2898e-03 & 2.80 & 1.4608e-03 & 1.97 \\
    4 & 4.4194e-02 & 163585 & 1.8953e-02 & 1.82 & 4.6687e-04 & 2.82 & 3.6503e-04 & 2.00 \\
    \bottomrule
  \end{tabular}%
  \caption{Convergence results for the HYPRE method using the Raviart--Thomas interpolator defined in Listing \ref{lst:raviart.thomas.interpolator} applied to the Kovasznay problem.}
  \label{tab:kovasznay:hypre}
\end{table}

Table \ref{tab:kovasznay:hypre} displays convergence results for the HYPRE method of \cite{Beirao-da-Veiga.Di-Pietro.ea:25} applied to the Kovasznay problem, which was already considered using the HHO method in Section \ref{sec:advanced.features:cartesian.product.spaces}.
These results are directly comparable with those in Table \ref{tab:kovasznay:hho} because they were obtained using the same sequence of refined meshes.
For $k = 0$, it is interesting to observe a transition from the convection-dominated to the diffusion-dominated order of convergence.

\subsection{Orthogonal complement polynomial spaces}

In Listing \ref{lst:potential}, the HHO potential is defined by enforcing conditions on the $L^2$-orthogonal decomposition $\text{\lstinline{Poly(k+1, scalar)}} = \text{\lstinline{ZeroAveragePoly(k+1, scalar)}} \oplus_{\perp_{L^2(T)}} \text{\lstinline{Poly(0, scalar)}}$.
More general $L^2$-orthogonal decompositions can be defined in the DSL, as shown in the following example, which presents the deflection reconstruction used in the HHO method for the Kirchhoff--Love problem (see \cite[Eqs. (3.4)--(3.5)]{Bonaldi.Di-Pietro.ea:18}):
\begin{lstlisting}
operator deflection_reconstruction : Uh(u) -> Poly(k+2, scalar) on element T {
  forall w in orthogonal complement of Poly(1, scalar) relative to Poly(k+2, scalar):
    int(T) grad(grad(deflection_reconstruction(u))) dot grad(grad(w)) =
    int(T) dof(u, T) * div(div(grad(grad(w))))
    + int(dT) dof(u, E, edge_gradient) dot (grad(grad(w)) dot normal)
    - int(dT) dof(u, E, edge_value) * (div(grad(grad(w))) dot normal)
  forall w in Poly(1, scalar):
    int(T) deflection_reconstruction(u) * w = int(T) dof(u, T) * w
}
\end{lstlisting}

\subsection{Mixed element and boundary forms}\label{sec:advanced.features:boundary.forms}

The DSL seamlessly handles the algebraic sum of element and boundary forms.
As an example, consider the following listing, which contains the modified forms for the HHO method discussed in Section \ref{sec:overview} with weakly enforced boundary conditions.
Notice, in passing, the syntax of \lstinline{penalty_weight}, regarded as a (constant-valued) function of the edge.
\begin{lstlisting}
parameter penalty = 10.0
  
function penalty_weight(E) -> scalar = penalty * pow(diameter(E), -1.0)

bilinear form poisson : Uh(trial u) times Uh(test v) {
  sum_elements(
    int(T) gradient_reconstruction(u) dot gradient_reconstruction(v)
    + pow(diameter(T), -2.0) * int(T) element_difference(u) * element_difference(v)
    + pow(diameter(T), -1.0) * int(dT) edge_difference(u) * edge_difference(v)
  )
  + sum_boundary_edges(
      - int(E) (gradient_reconstruction(u) dot normal) * dof(v, E)
      - int(E) dof(u, E) * (gradient_reconstruction(v) dot normal)
      + penalty_weight(E) * int(E) dof(u, E) * dof(v, E)
    )
}

linear form load : Uh(test v) {
  sum_elements(int(T) source * dof(v, T))
  + sum_boundary_edges(
      - int(E) exact_solution * (gradient_reconstruction(v) dot normal)
      + penalty_weight(E) * int(E) exact_solution * dof(v, E)
    )
}
\end{lstlisting}


\section{Conclusion}\label{sec:conclusion}

\polyrust provides a unified environment for expressing and testing polytopal discretizations in a syntax that stays close to their fully discrete mathematical formulation.
Its external DSL covers the complete path from discrete spaces and interpolants to local reconstructions, forms, boundary conditions, algebraic problems, error measures, and field export.
The same language accommodates a variety of methods, including, for example, HHO, DDR, and Virtual Element methods, as well as linear and nonlinear problems set on atomic or Cartesian product spaces.
Combined with an efficient and memory-safe \Rust implementation, automated exactness checks, parallel assembly, and convergence-study tools, this makes \polyrust particularly well suited to the rapid and reliable prototyping of new polytopal methods.
Its principal asset is therefore the ability to turn complex operator-based formulations into concise, reusable, and directly executable descriptions.


\section*{Acknowledgements}

This project has received funding from the European Research Council (ERC) under the European Union’s Horizon Europe research and innovation programme (Grant agreement No. 101115663).

Funded by the European Union (ERC Synergy, NEMESIS, project number 101115663).
Views and opinions expressed are however those of the authors only and do not necessarily reflect those of the European Union or the European Research Council Executive Agency. Neither the European Union nor the granting authority can be held responsible for them.

The author also wishes to thank Louis Lam\'erand for providing the refined sequence of polyhedral meshes used in numerical tests.


\section*{Data availability}

This article does not make use of data.
All numerical experiments can be replicated using the \polyrust code.
Instructions are available at \url{https://polyrust.org}.

\section*{Ethics and consent}

Ethical approval and consent were not required.

\section*{Competing interests}

The author has no competing interests to declare.

\section*{Software and code availability}
  
The \polyrust code is available through the public GIT repository
\url{git@plmlab.math.cnrs.fr:dpietro/poly_rust-releases.git}
Packages for Mac OS and Linux are available at \url{https://polyrust.org}.
Version 0.13.1 is considered for the present work.


\printbibliography

@Article{         Alns.Logg.ea:14,
  author        = {Aln\ae{}s, Martin S. and Logg, Anders and {\O}lgaard,
                  Kristian B. and Rognes, Marie E. and Wells, Garth N.},
  title         = {Unified Form Language: A domain-specific language for weak
                  formulations of partial differential equations},
  journal       = {ACM Trans. Math. Softw.},
  volume        = {40},
  number        = {2},
  pages         = {9:1--9:37},
  year          = {2014},
  doi           = {10.1145/2566630}
}

@InCollection{    Antonietti.Cangiani.ea:16,
  author        = {Antonietti, P. F. and Cangiani, A. and Collis, J. and
                  Dong, Z. and Georgoulis, E. H. and Giani, S. and Houston,
                  P.},
  title         = {Review of discontinuous {G}alerkin finite element methods
                  for partial differential equations on complicated domains},
  booktitle     = {Building bridges: connections and challenges in modern
                  approaches to numerical partial differential equations},
  series        = {Lect. Notes Comput. Sci. Eng.},
  volume        = {114},
  pages         = {279--308},
  publisher     = {Springer, [Cham]},
  year          = {2016}
}

@Article{         Antonietti.Giani.ea:13,
  author        = {Antonietti, P. F. and Giani, S. and Houston, P.},
  title         = {$hp$-version composite discontinuous {G}alerkin methods
                  for elliptic problems on complicated domains},
  journal       = {SIAM J. Sci. Comput.},
  volume        = {35},
  number        = {3},
  pages         = {A1417--A1439},
  year          = {2013},
  doi           = {10.1137/120877246}
}

@Article{         Badia.Verdugo:20,
  author        = {Badia, Santiago and Verdugo, Francesc},
  title         = {Gridap: An extensible finite element toolbox in {J}ulia},
  journal       = {J. Open Source Softw.},
  volume        = {5},
  number        = {52},
  pages         = {2520},
  year          = {2020},
  doi           = {10.21105/joss.02520}
}

@Article{         Bangerth.Hartmann.ea:07,
  author        = {Bangerth, Wolfgang and Hartmann, Ralf and Kanschat,
                  Guido},
  title         = {deal.{II}---A general-purpose object-oriented finite
                  element library},
  journal       = {ACM Trans. Math. Softw.},
  volume        = {33},
  number        = {4},
  pages         = {24:1--24:27},
  year          = {2007},
  doi           = {10.1145/1268776.1268779}
}

@Article{         Bassi.Botti.ea:12,
  author        = {Bassi, F. and Botti, L. and Colombo, A. and Di Pietro, D.
                  A. and Tesini, P.},
  title         = {On the flexibility of agglomeration based physical space
                  discontinuous {Galerkin} discretizations},
  journal       = {J. Comput. Phys.},
  volume        = {231},
  number        = {1},
  pages         = {45--65},
  year          = {2012},
  doi           = {10.1016/j.jcp.2011.08.018}
}

@Misc{            Beaude.Lemaire:21,
  title         = {{ParaSkel++: a C++ platform for the high-performance,
                  arbitrary-order, 2/3D numerical approximation of PDEs on
                  general polytopal meshes using skeletal Galerkin methods}},
  author        = {Beaude, Laurence and Lemaire, Simon},
  url           = {https://hal.science/hal-03517921},
  year          = {2021}
}

@Article{         Beirao-da-Veiga.Brezzi.ea:13,
  author        = {{Beir\~{a}o da Veiga}, L. and Brezzi, F. and Cangiani, A.
                  and Manzini, G. and Marini, L. D. and Russo, A.},
  title         = {Basic principles of virtual element methods},
  journal       = {Math. Models Methods Appl. Sci.},
  number        = {23},
  volume        = {199},
  year          = {2013},
  pages         = {199--214},
  doi           = {10.1142/S0218202512500492}
}

@Article{         Beirao-da-Veiga.Brezzi.ea:14,
  title         = {The hitchhiker's guide to the virtual element method},
  author        = {{Beir{\~a}o da Veiga}, L. and Brezzi, F. and Marini, L.D.
                  and Russo, A.},
  journal       = {Math. Models Methods Appl. Sci.},
  volume        = {24},
  number        = {08},
  pages         = {1541--1573},
  year          = {2014},
  doi           = {10.1142/S021820251440003X}
}

@Article{         Beirao-da-Veiga.Brezzi.ea:23,
  title         = {The virtual element method},
  author        = {{Beir\~ao da Veiga}, Louren\c{c}o and Brezzi, Franco and
                  Marini, L Donatella and Russo, Alessandro},
  journal       = {Acta Numerica},
  volume        = {32},
  pages         = {123--202},
  year          = {2023},
  doi           = {10.1017/S0962492922000095}
}

@Article{         Beirao-da-Veiga.Di-Pietro.ea:25,
  author        = {{Beir\~{a}o da Veiga}, L. and Di Pietro, D. A. and
                  Droniou, J. and Haile, K. B. and Radley, T. J.},
  title         = {A {Reynolds}-semi-robust method with hybrid velocity and
                  pressure for the unsteady incompressible {Navier--Stokes}
                  equations},
  journal       = {SIAM J. Numer. Anal.},
  year          = {2025},
  volume        = {63},
  number        = {6},
  pages         = {2317--2342},
  doi           = {10.1137/25M1736104}
}

@Misc{            Beirao-da-Veiga.Di-Pietro.ea:26,
  title         = {Key challenges and bridges among convergence analysis
                  techniques for polytopal methods},
  author        = {{Beir\~{a}o da Veiga}, L. and Di Pietro, D. A. and
                  Droniou, J.},
  year          = {2026},
  month         = {5},
  eprint        = {2605.23405},
  archiveprefix = {arXiv},
  primaryclass  = {math.NA}
}

@Article{         Berrone.Borio.ea:26,
  author        = {Berrone, Stefano and Borio, Andrea and Teora, Gioana and
                  Vicini, Fabio},
  title         = {{POLYDIM}: A {C++} library for {POLY}topal
                  {DI}scretization {M}ethods},
  journal       = {Comput. Phys. Commun.},
  volume        = {320},
  pages         = {109937},
  year          = {2026},
  doi           = {10.1016/j.cpc.2025.109937}
}

@Article{         Bonaldi.Di-Pietro.ea:18,
  author        = {Bonaldi, F. and Di Pietro, D. A. and Geymonat, G. and
                  Krasucki, F.},
  title         = {A {Hybrid High-Order} method for {Kirchhoff--Love} plate
                  bending problems},
  journal       = {ESAIM: Math. Model Numer. Anal.},
  volume        = {52},
  number        = {2},
  pages         = {393--421},
  year          = {2018},
  doi           = {10.1051/m2an/2017065}
}

@Article{         Bonaldi.Di-Pietro.ea:25,
  title         = {An exterior calculus framework for polytopal methods},
  author        = {Bonaldi, F. and Di Pietro, D. A. and Droniou, J. and Hu,
                  K.},
  journal       = {J. Eur. Math. Soc.},
  year          = {2025},
  doi           = {10.4171/JEMS/1602},
  note          = {Published online}
}

@Article{         Botti.Botti.ea:26,
  author        = {Botti, L. and Botti, M. and Di Pietro, D. A. and Massa, F.
                  C.},
  title         = {Stability, convergence, and pressure-robustness of
                  numerical schemes for incompressible flows with hybrid
                  velocity and pressure},
  year          = {2026},
  journal       = {Math. Comp.},
  volume        = {95},
  number        = {357},
  pages         = {1--28},
  doi           = {10.1090/mcom/4049}
}

@Article{         Botti.Di-Pietro.ea:19*1,
  author        = {Botti, L. and Di Pietro, D. A. and Droniou, J.},
  title         = {A {Hybrid High-Order} method for the incompressible
                  {Navier--Stokes} equations based on {Temam}'s device},
  year          = {2019},
  volume        = {376},
  pages         = {786--816},
  doi           = {10.1016/j.jcp.2018.10.014},
  journal       = {J. Comput. Phys.}
}

@Article{         Chin.Lasserre.ea:15,
  author        = {Chin, E. B. and Lasserre, J. B. and Sukumar, N.},
  title         = {Numerical integration of homogeneous functions on convex
                  and nonconvex polygons and polyhedra},
  journal       = {Comput. Mech.},
  volume        = {56},
  year          = {2015},
  number        = {6},
  pages         = {967--981},
  doi           = {10.1007/s00466-015-1213-7}
}

@Article{         Cicuttin.Di-Pietro.ea:18*1,
  author        = {Cicuttin, M. and Di Pietro, D. A. and Ern, A.},
  title         = {Implementation of {Discontinuous Skeletal} methods on
                  arbitrary-dimensional, polytopal meshes using generic
                  programming},
  journal       = {J. Comput. Appl. Math.},
  volume        = {344},
  pages         = {852--874},
  year          = {2018},
  doi           = {10.1016/j.cam.2017.09.017}
}

@InBook{          Cicuttin.Ern.ea:21,
  author        = {Cicuttin, Matteo and Ern, Alexandre and Pignet, Nicolas},
  title         = {Implementation Aspects},
  booktitle     = {Hybrid High-Order Methods: A Primer with Applications to
                  Solid Mechanics},
  year          = {2021},
  publisher     = {Springer International Publishing},
  address       = {Cham},
  pages         = {109--127},
  doi           = {10.1007/978-3-030-81477-9_8}
}

@Article{         Dassi:26,
  author        = {Dassi, Franco},
  title         = {{Vem++}, a {C++} library to handle and play with the
                  virtual element method},
  journal       = {Numer. Algorithms},
  volume        = {101},
  pages         = {1633--1675},
  year          = {2026},
  doi           = {10.1007/s11075-025-02059-z}
}

@Article{         Davis:04,
  author        = {Davis, Timothy A.},
  title         = {Algorithm 832: {UMFPACK} V4.3---an unsymmetric-pattern
                  multifrontal method},
  journal       = {ACM Trans. Math. Softw.},
  volume        = {30},
  number        = {2},
  pages         = {196--199},
  year          = {2004},
  doi           = {10.1145/992200.992206}
}

@Article{         Di-Pietro.Droniou.ea:15,
  author        = {Di Pietro, D. A. and J. Droniou and Ern, A.},
  title         = {A discontinuous-skeletal method for
                  advection-diffusion-reaction on general meshes},
  year          = {2015},
  journal       = {SIAM J. Numer. Anal.},
  volume        = {53},
  number        = {5},
  pages         = {2135--2157},
  doi           = {10.1137/140993971}
}

@Article{         Di-Pietro.Droniou.ea:20,
  author        = {Di Pietro, D. A. and Droniou, J. and Rapetti, F.},
  title         = {Fully discrete polynomial {de Rham} sequences of arbitrary
                  degree on polygons and polyhedra},
  journal       = {Math. Models Methods Appl. Sci.},
  year          = {2020},
  volume        = {30},
  number        = {9},
  pages         = {1809-1855},
  doi           = {10.1142/S0218202520500372}
}

@Misc{            Di-Pietro.Droniou.ea:25,
  author        = {Di Pietro, D. A. and Droniou, J. and Hu, K. and Leroy,
                  A.},
  title         = {Analytical properties of polygonal {Stokes} and {BGG}
                  Hessian complexes with application to {Kirchhoff--Love}
                  plates},
  year          = {2025},
  month         = {7},
  eprint        = {2507.17333},
  archiveprefix = {arXiv},
  primaryclass  = {math.NA}
}

@Book{            Di-Pietro.Droniou:20,
  author        = {Di Pietro, D. A. and Droniou, J.},
  title         = {The {Hybrid High-Order} method for polytopal meshes},
  subtitle      = {Design, analysis, and applications},
  publisher     = {Springer International Publishing},
  year          = {2020},
  series        = {Modeling, Simulation and Application},
  number        = {19},
  doi           = {10.1007/978-3-030-37203-3}
}

@Article{         Di-Pietro.Droniou:23*2,
  author        = {Di Pietro, D. A. and Droniou, J.},
  title         = {An arbitrary-order discrete {de Rham} complex on
                  polyhedral meshes: Exactness, {Poincar\'e} inequalities,
                  and consistency},
  journal       = {Found. Comput. Math.},
  year          = {2023},
  volume        = {23},
  pages         = {85--164},
  doi           = {10.1007/s10208-021-09542-8}
}

@Article{         Di-Pietro.Ern.ea:14,
  author        = {Di Pietro, D. A. and Ern, A. and Lemaire, S.},
  title         = {An arbitrary-order and compact-stencil discretization of
                  diffusion on general meshes based on local reconstruction
                  operators},
  journal       = {Comput. Meth. Appl. Math.},
  volume        = {14},
  number        = {4},
  pages         = {461--472},
  year          = {2014},
  doi           = {10.1515/cmam-2014-0018}
}

@Article{         Di-Pietro.Ern.ea:16*1,
  author        = {Di Pietro, D. A. and Ern, A. and Linke, A. and Schieweck,
                  F.},
  title         = {A discontinuous skeletal method for the
                  viscosity-dependent {Stokes} problem},
  journal       = {Comput. Meth. Appl. Mech. Engrg.},
  year          = {2016},
  volume        = {306},
  pages         = {175--195},
  doi           = {10.1016/j.cma.2016.03.033}
}

@Article{         Di-Pietro.Ern:10,
  author        = {Di Pietro, D. A. and Ern, A.},
  title         = {Discrete functional analysis tools for discontinuous
                  {G}alerkin methods with application to the incompressible
                  {N}avier-{S}tokes equations},
  journal       = {Math. Comp.},
  volume        = {79},
  year          = {2010},
  number        = {271},
  pages         = {1303--1330},
  doi           = {10.1090/S0025-5718-10-02333-1}
}

@Article{         Di-Pietro.Ern:15,
  author        = {Di Pietro, D. A. and Ern, A.},
  title         = {A hybrid high-order locking-free method for linear
                  elasticity on general meshes},
  journal       = {Comput. Meth. Appl. Mech. Engrg.},
  year          = {2015},
  volume        = {283},
  pages         = {1--21},
  doi           = {10.1016/j.cma.2014.09.009}
}

@Article{         Droniou.Yemm:22,
  author        = {Droniou, J\'{e}r\^{o}me and Yemm, Liam},
  title         = {Robust hybrid high-order method on polytopal meshes with
                  small faces},
  journal       = {Comput. Methods Appl. Math.},
  volume        = {22},
  year          = {2022},
  number        = {1},
  pages         = {47--71},
  doi           = {10.1515/cmam-2021-0018}
}

@Article{         Hecht:12,
  author        = {Hecht, Fr\'{e}d\'{e}ric},
  title         = {New development in {FreeFem++}},
  journal       = {J. Numer. Math.},
  volume        = {20},
  number        = {3--4},
  pages         = {251--265},
  year          = {2012},
  doi           = {10.1515/jnum-2012-0013}
}

@Article{         Kovasznay:48,
  author        = {Kovasznay, L. I. G.},
  title         = {Laminar flow behind a two-dimensional grid},
  journal       = {Math. Proc. Camb. Philos. Soc.},
  year          = {1948},
  volume        = {44},
  number        = {1},
  pages         = {58--62},
  doi           = {10.1017/S0305004100023999},
  publisher     = {Cambridge University Press}
}

@Article{         Prudhomme.Chabannes.ea:12,
  author        = {Prud'homme, Christophe and Chabannes, Vincent and Doyeux,
                  Vincent and Ismail, Mourad and Samake, Abdoulaye and Pena,
                  Gon\c{c}alo},
  title         = {{Feel++}: A computational framework for {Galerkin} methods
                  and advanced numerical methods},
  journal       = {ESAIM Proc.},
  volume        = {38},
  pages         = {429--455},
  year          = {2012},
  doi           = {10.1051/proc/201238024}
}

\end{document}